\documentclass[draft,12pt,leqno]{amsart}

\usepackage{amsmath,amssymb,amscd,amsthm}
\usepackage{latexsym}
\usepackage[dvips]{graphics,epsfig}

\newtheorem{theorem}{Theorem}
\newtheorem{lemma}{Lemma}
\newtheorem{proposition}{Proposition}

\newtheorem{definition}{Definition}
\newtheorem{corollary}{Corollary}

\newtheorem{claim}{Claim}

\newcommand{\qq}{\quad\quad}

\newcommand{\norm}[2]{{\left\| #1 \right\|}_{#2}}
\newcommand{\f}[2]{\frac{#1}{#2}}
\newcommand{\dpr}[2]{\langle #1,#2 \rangle}

\newcommand{\al}{\alpha}

\newcommand{\de}{\delta}

\newcommand{\ve}{\varepsilon}

\newcommand{\La}{\Lambda}
\newcommand{\si}{\sigma}

\newcommand{\vp}{\varphi}

\newcommand{\rn}{{\mathbf R}^n}

\newcommand{\rone}{\mathbf R^1}

\newcommand{\rthree}{\mathbf R^3}

\newcommand{\ca}{\mathcal A}
\newcommand{\cs}{\mathcal S}

\newcommand{\cz}{\mathcal Z}

\newcommand{\cc}{\mathcal C}

\newcommand{\intl}{\int\limits}

\newcommand{\liml}{\lim\limits}
\newcommand{\suml}{\sum\limits}

\newcommand{\supl}{\sup\limits}

\newcommand{\p}{\partial}

\newcommand{\ed}{(1-\p_x^2)^{-1}}

\newcommand{\beq}{\begin{equation}}
\newcommand{\eeq}{\end{equation}}
\newcommand{\beqna}{\begin{eqnarray*}}
\newcommand{\eeqna}{\end{eqnarray*}}
\newcommand{\beqn}{\begin{equation*}}
\newcommand{\eeqn}{\end{equation*}}
\newcommand{\bp}{\begin{proof}}
\newcommand{\ep}{\end{proof}}
\newcommand{\bprop}{\begin{proposition}}
\newcommand{\eprop}{\end{proposition}}
\newcommand{\bt}{\begin{theorem}}
\newcommand{\et}{\end{theorem}}
\newcommand{\bex}{\begin{Example}}
\newcommand{\eex}{\end{Example}}
\newcommand{\bc}{\begin{corollary}}
\newcommand{\ec}{\end{corollary}}
\newcommand{\bcl}{\begin{claim}}
\newcommand{\ecl}{\end{claim}}
\newcommand{\bl}{\begin{lemma}}
\newcommand{\el}{\end{lemma}}

\begin{document}

\title
[Attractors for the viscous Camassa-Holm equation]
{Attractors for the viscous Camassa-Holm equation}

\author{Milena Stanislavova} 
\author{Atanas Stefanov}

\address{Milena Stanislavova\\
Department of Mathematics \\
University of Kansas\\
1460 Jayhawk Blvd\\ Lawrence, KS 66045--7523}
\email{stanis@math.ku.edu}
\address{Atanas Stefanov\\
Department of Mathematics \\
University of Kansas\\
1460 Jayhawk Blvd\\ Lawrence, KS 66045--7523}

\email{stefanov@math.ku.edu}

\thanks{
First author supported in part by the NSF under grant \# EPS-0236913 and 
NSF-DMS 0508184. 
Second author supported in part by  
NSF-DMS 0300511.}
\date{\today}



\begin{abstract}
We consider the viscous Camassa-Holm equation subject to an external force, 
where the viscosity term is   
given by second order differential operator in  divergence form. We show that under some mild assumptions on the viscosity term, one has global well-posedness both in 
the periodic case and the case of the whole line. In the periodic case, 
we show  the existence of  global attractors in the energy space $H^1$, 
provided the external force is in the class $L^2(I)$.  Moreover, we establish 
an asymptotic smoothing effect, which states that the elements of the attractor are 
in fact in the smoother Besov space $B^2_{2,\infty}(I)$.  
Identical  results (after  adding an appropriate linear damping  term) 
are obtained in the case of the whole line. 
\end{abstract}

\maketitle

\section{Introduction}
The failure of weakly nonlinear dispersive equations, such as 
the celebrated Korteweg-de Vries equation, 
to model interesting physical phenomena like wave breaking, 
existence of peaked waves etc., was a motivation for transition to full nonlinearity in 
the search for alternative models for nonlinear dispersive waves (\cite{Whitham}). 
The first step in this direction was the derivation of the Green-Naghdi system of equations (see \cite{GN}), which is a Hamiltonian system that models fluid flows in thin domains. Writing the Green-Naghdi equations in Hamiltonian form and using asymptotic expansion which keeps the Hamiltonian structure, Camassa and Holm (\cite{Camassa}) derived the Camassa-Holm equation in 1993. They obtained the strongly nonlinear equation
$$
u_t-\frac{1}{4} u_{xxt}+\frac{3}{2} (u^2)_x-\frac{1}{8} (u_x^2)_x-\frac{1}{4} (u u_{xx})_x=0,
$$
which was also found independently by Dai (\cite{Dai}) as a model for nonlinear waves in 
cylindrical hyper elastic rods and had been originally obtained by Fokas and Fuchsteiner 
(\cite{FF}) as an example of bi-Hamiltonian equation. The equation possesses 
a Lax pair and is  completely integrable in terms of the inverse scattering transform, 
see \cite{Cam1},\cite{Camassa}.
For recent and 
extensive treatments of the case of solutions on the real line, decaying at infinity, 
we refer to \cite{con5}, \cite{con1}, 
\cite{kaup}, while for the periodic case, one should consult \cite{con2}, \cite{con4} .

 A dictinct feature of the Camassa-Holm equations is that it 
 exhibits orbitally stable soliton solutions, which are weak solutions in the shape of a 
peaked waves, \cite{Co2}, \cite{Co1}, see also \cite{lenells} 
for the most  complete description of traveling waves available as of 
this writing. 
 Camassa and Holm, \cite{Camassa} have found that two 
solitary waves keep their shape and size after interaction while the ultimate 
position of each wave is affected only with a phase shift by the nonlinear interaction, 
see also \cite{Beals},
\cite{Co3}. Finally we mention the presence of breaking waves for this equation 
(\cite{Camassa},  \cite{Co},\cite{McKean}, \cite{con3}), as well as the occurrence of global solutions 
(\cite{Co}, \cite{Co2}, \cite{Co3},\cite{Co4}).

Our main object of investigation will be 
 the initial value problem for 
the Camassa-Holm equation, which takes the form 
\begin{equation}
\label{eq:1}
(CH) \ \ \left|\begin{array}{l}
u_t-u_{txx}=2u_x u_{xx}+u u_{xxx}-3 u u_x\\
u(x,0)=u_0(x)
\end{array}\right.
\end{equation}
By reorganizing the terms, one sees that this is equivalent to
\begin{equation}
\label{eq:10}
u_t+\f{1}{2}\p_x(u^2)+\p_x\ed [u_x^2/2+u^2]=0,
\end{equation}
where the Helmholtz operator $\ed$ is standardly defined in Section \ref{sec:helm}. 
Denote here and for the rest of the paper the nonlinearity of \eqref{eq:10} 
 $F(u,u_x)= \f{1}{2}\p_x(u^2)+\p_x\ed [u_x^2/2+u^2]$. 

The  viscous Camassa-Holm equation in one and more 
dimensions\footnote{These equations are also known as 
 Navier Stokes $\al$ models.} was studied extensively in the recent years.  
This was done in parallel with  the non viscous one, so we refer to the papers, 
quoted above. 
In \cite{Stanislavova}, we have shown in particular that for 
\begin{equation}
\label{eq:11}
u_t+\f{1}{2}\p_x(u^2)+\p_x\ed [u_x^2/2+u^2]=\ve \p_x^2 u,
\end{equation}
one has global and unique solution in the energy class $H^1(\rone)$. 

In \cite{Helge}, the authors have taken a more general type of 
viscosity and forcing terms. They have shown (among other 
things) global well-posedness for the equation 
\begin{equation}
\label{eq:2}
u_t+\f{1}{2}\p_x(u^2)+\p_x\ed [u_x^2/2+u^2]=\p_x(a(t,x)\p_x  u)+g(t,x,u),
\end{equation}
with initial data $u_0\in H^2(\rone)$. 
Here $a$ is bounded, positive and bounded away from 
zero, with number of additional technical assumptions on $a$, $g$. 

In this article, we shall consider  similar type of 
 viscosity terms $\p_x(a \p_x u)$, 
which is motivated by recent works in conservation laws and 
which seem to better model the underlying physical situations. 
We will  however stick to the case of {\it time independent} $a=a(x)$, 
although our arguments work in the time dependent case as well, 
subject to some minor modifications. This is done to reduce the unnecessary 
technicalities and 
it is also dictated by our interest in  the dynamical system   
(rather than the cocycle) properties of \eqref{eq:2}.  

It is also our goal to consider the question for global 
well-posedness of \eqref{eq:2} both on the whole line 
$\rone$ and on any finite interval. As we shall see, 
 the methods that we employ in the two cases  are slightly different, 
but not conceptually so. The main difficulty for the case of $\rone$ 
 as usual 
will be the non compactness of the embedding 
$H^1(\rone)\hookrightarrow L^2(\rone)$. 

Let us take a moment to explain our results. First, 
under standard  assumptions\footnote{In fact, for the existence theorem, 
the smoothness assumptions on $a$ that we work with 
are considerably less restrictive than those imposed by \cite{Helge}. 
Moreover, in the proof of well-posedness, 
it will suffice to assume only $a\in C^1(I)$.} 
on $a$ and $g$, 
{\it we show that the dynamical system \eqref{eq:2} has an 
unique global solution}, whenever $u_0\in H^1(\rone)$ or $u_0\in H^1(0, 1)$ respectively. 

For the case of finite interval, 
 we are  able to show {\it the existence 
of global attractor}. 
This is done under a smallness assumption on the Lipschitz norm of 
$a$. 

In addition,  the attractor (which is initially a subset of $H^1(0, 1)$) 
turns out to be a subset of the  
smoother space  $H^{2-\si}(0, 1)$,  {\it that is the 
semigroup associated with \eqref{eq:2}  exhibits asymptotic smoothing effect}. More precisely, 
we show that for every $\si>0$, the attractor is a bounded subset of $H^{2-\si}(0,1)$. 

For the case of \eqref{eq:2} considered as a integro differential equation 
on the whole line 
 $\rone$, the existence of an attractor is not clear, although 
we have not explicitly found a counterexample. 
The main difficulty is that nothing seems to prevent a 
low-frequency buildup, which may cause an unrestricted growth of 
$\|u(t, \cdot)\|_{H^1}$. That is, we expect that for a wide class of 
 initial data $u_0$ and right hand side $g$, 
$\limsup_{t\to \infty} \|u(t, \cdot)\|_{H^1(\rone)}=\infty$. This 
clearly would prevent the existence of an attractor. 

On the other hand, 
if one adds an 
additional damping term (which is 
actually a relevant physical model, considered in two dimensions 
by Ilyin and Titi, \cite{Titi5}), one can show 
the existence of an attractor and  boundedness in $H^{2-}$
 in the case of the whole line as well. The 
discussion on that is in Section \ref{sec:conclusions}. \\
Now and throughout the paper, we will require that the operator 
$A u=-\p_x(a(x) \p_x u)$ be  coercive. 
That is, assume that   $a(x)$ is $C^2$ real valued, so that  
 for some fixed $\ve>0$
\begin{equation}
\label{eq:c1}
\ve<a(x) <1/\ve 
\end{equation}
Note that under these assumptions, we can define  
the (unbounded) operator $A$ as a Friedrich's extension of 
the unbounded operator defined by  the quadratic 
 form 
$$
q(u,u) = \int_I a(x) |u'(x)|^2 dx=:\dpr{u}{Au},
$$
with domain $\dot{H}^1(I)$ with the natural boundary conditions and 
where $I=\rone$ or $I=(0,1)$.  That is, 
we impose the  boundary condition $u(0)=u(1)$ 
 in the periodic case and $\lim_{|x|\to \infty} u(x)=0$ in the case of the whole line. 
In particular, $A$ is positive and 
self-adjoint operator and $-A$ generates a strongly continuous semigroup. \\
Our first theorem is a well-posedness type result. 
\begin{theorem}
\label{theo:1} For the viscous Camassa-Holm equation \eqref{eq:2},  
assume that $a=a(x)$ satisfies\footnote{For the well-posedness result, 
it is enough to assume only 
that $a\in C^1(0,1)$.} 
\eqref{eq:c1} and $g\in L^\infty_t L^2_x (I)$, where either 
$I=\rone$ or $I=(0,1)$. Then for every initial data $u_0\in H^1(I)$, 
there is an unique global classical solution $u$ to \eqref{eq:2}. 
More specifically, $u\in C([0,\infty), H^1(I))$ and for every 
$0<T_1<T_2<\infty$, $u\in C^2([T_1, T_2], I)$. 
\end{theorem}
Our next result concerns the existence of global attractors 
for \eqref{eq:2} in the case of finite 
interval\footnote{As we have mentioned already, 
global attractors may not exist in the case $I=\rone$.}  
$I=(0,1)$.  For technical reasons, we need to impose a  smallness condition 
 $\norm{a'}{L^\infty}<<\ve$. We do not know whether such a condition is 
necessary or not, but it is possible that unless such a condition hold, 
one gets  unbounded orbits for some 
sets of initial data, thus rendering the statements regarding the 
existence of attractors false.  
\begin{theorem}
\label{theo:2} Assume that $a$ satisfies 
\eqref{eq:c1} and $\norm{a'}{L^\infty}\leq \de \ve$ for some sufficiently small $\de$. 
Let  $g=g(x)\in L^2(0,1)$ has mean value zero, $\int_0^1 g(x) dx=0$. 
Then, the viscous Camassa-Holm equation \eqref{eq:2} 
has a global attractor, when considered as a dynamical system over a 
finite interval $I=(0,1)$ with initial data in 
$H^1_0(0,1)=H^1\cap 
\{f: \int_0^1  f(x) dx=0\}$. 
\end{theorem}
\noindent {\bf Remark:} The mean value zero condition imposed upon the 
forcing term $g$ is necessary 
for the existence of a global attractor and is in fact  necessary merely  for uniform 
boundedness of the orbits. 

Indeed, an elementary computation shows that $
\p_t \int_0^1   u (t,x) dx= \int_0^1 g(x) dx$, 
whence $\int_0^1 u(t,x) dx = \int_0^1 u_0(x) dx+ (\int g(x) dx)  t$, which is not bounded as 
$t\to \infty$, unless   $\int g(x) dx=0$. 

Our next theorem addresses 
 precisely the asymptotic smoothing effect of the corresponding dynamics.  
\begin{theorem}
\label{theo:3}
The attractor $\ca$ constructed in Theorem \ref{theo:2} is contained in 
$\cap_{\si>0} H^{2-\si}(0,1)$. Moreover, for all $\si>0$, we have the estimate 
$$
\supl_{f \in \ca}  \norm{f}{H^{2-\si}(0,1)}\leq C_\si \norm{g}{L^2}
$$
That is, the attractor is a bounded subset in $H^{2-\si}$ 
with bounds depending only on  the constants in the problem ($\ve, \de, \si$) and $\norm{g}{L^2}$. 

In fact, more generally,   
$\ca$ is a bounded subset of  $ B^{2}_{2, \infty}$, with the 
corresponding estimate 
\begin{equation}
\label{eq:smooth}
\supl_{f\in \ca} \supl_k 2^{2k} \norm{P_{2^k} f}{L^2}= \supl_{f\in \ca} \supl_k  2^{2k} 
\left(\suml_{n=2^{k-1}}^{2^{k+1}} |\hat{f}(n)|^2\right)^{1/2} 
\leq C\norm{g}{L^2},
\end{equation}
\end{theorem}
 We record that in the case of constant viscosity (i.e. $a=const>0$), 
all the conditions in Theorem \ref{theo:2} and Theorem \ref{theo:3} 
are satisfied. 

For the case of the whole line,  consider the 
Camassa-Holm equation with an additional damping factor, 
as considered in two dimensions by Ilyin-Titi, \cite{Titi5}. Namely, let 
$\mu>0$ and consider 
\begin{equation}
\label{e:1}
u_t +\mu u + \f{1}{2}\p_x(u^2)+\p_x\ed [u_x^2/2+u^2] =\p_x(a(x) u_x)+g(x)
\end{equation}
with initial data $u(0,x)=f$. We have the following 
\begin{theorem}
\label{theo:9}
Assume that $a$ satisfies \eqref{eq:c1} and either 
$\norm{a'}{L^\infty}<\de\ve$ for some sufficiently 
small $\de$ {\it or} $a''(x)\leq 2 a(x)$. 
Then the equation \eqref{e:1} is globally well-posed in $H^1(\rone)$. 
It also has a global attractor $\ca$ and the semigroup has the smoothing property: 
$\ca$ is a bounded subset of $B^2_{2, \infty}$. 
More precisely, 
$$
\supl_{f\in \ca}\sup_k 2^{2k} \norm{P_{2^k} f}{L^2}\leq C\norm{g}{L^2}.
$$ 
\end{theorem}
\noindent {\bf Remark}
\begin{itemize}
\item If  $a=const>0$, all the 
conditions in Theorem \ref{theo:9} are met 
and the results hold. 
\item In contrast with Theorem \ref{theo:2}, note that we can impose 
the structural condition $a''(x)\leq 2 a(x)$, instead of the 
smallness of $\norm{a'}{L^\infty}$.  
\end{itemize}
\noindent 
The paper is organized as follows. In Section \ref{sec:prelim}, we collect some useful facts
 from Fourier analysis and the theory of attractors both in finite and infinite domain setting. 
 In Section \ref{sec:90}, we first show a local  well-posedness of the Cauchy 
problem for the viscous Camassa-Holm equation, by using some elementary $C_0$ semigroup 
 properties of the semigroup generated by $A=-\p_x(a(\cdot) \p_x \cdot)$. This is 
done by a contraction map principle and yields valid solution only for short time. 
 We then derive\footnote{see Section \ref{sec:smoothness}} 
  some additional $H^2$ smoothness estimates in order to exploit the underlying $H^1$ 
  conservation law. 
  
  In Section \ref{sec:global}, we show that $H^1$ {\it a priori} estimates hold on 
{\it any time interval} $(0,T)$ and thus global well-posedness is established. 
  
  In Section \ref{sec:attractor_1}, we establish the existence of global 
attractors in the case of finite interval.  This is done by  verifying the 
point dissipativeness and the uniform boundedness of the dynamics.  
The uniform\footnote{Here uniform means uniformity with respect to a given 
bounded sequence of  initial data.} vanishing of the high frequency mass of 
the solutions, which is needed for the existence of attractors   
is addressed  in Section \ref{sec:vanish}.  Incidentally, 
  one obtains the smoothing estimate \eqref{eq:smooth}. 
  
  Finally, in Section \ref{sec:conclusions}, we prove Theorem \ref{theo:9}. The methods here are quite similar to the ones used in the final interval case. \\
\noindent   
{\bf Acknowledgement:} We are grateful to our colleague Bixiang Wang for 
numerous discussions regarding the abstract  theory of  attractors and their properties. 
 
\section{Preliminaries} 
\label{sec:prelim}
In this section, we collect some useful (generally well-known) 
 facts. 
We start with the definition of the Fourier transform 
in the whole space and in the periodic setting. 
\subsection{The Fourier transform and the Helmholtz operator $\ed$} 
\label{sec:helm}
The Fourier transform on $\rone$  is (initially) defined on the 
functions in the Schwartz class $\cs$ by 
$$
\hat{f}(\xi)=\intl_{\rone} f(x) e^{-2\pi i x\xi} dx.
$$
We record the inverse Fourier transform 
$$
f(x)=\intl_{\rone} \hat{f}(\xi) e^{2\pi i  x\xi } d\xi,
$$
and the Plancherel's identity is  $\norm{f}{L^2}= 
\|\hat{f}\|_{L^2}$ for all functions $f\in L^2$.  

On the interval  $[0,1]$, we may introduce the  Fouier transform 
$L^2([0,1])\to l^2(\cz)$, 
by setting $f\to \{a_k\}_{k\in\cz}$, where 
$$
a_k= \intl_0^1  f(x) e^{-2\pi i k x} dx.
$$
The inverse Fourier transform in that case is the familiar Fourier expansion 
$$
f(x)=\suml_{k \in \cz}  a_k e^{2\pi i k x}.
$$
an the Plancherel's identity is 
$\norm{f}{L^2([0,L])}=\norm{\{a_k\}}{l^2}$. 
Note that here and for the rest of the 
paper $L^2([0,1])$ is the space of square integrable functions with period one.

The Helmholtz operator is the inverse of the operator $(1-\p_x^2)$ or $\ed$. 
This is well-defined on both $L^2(\rone)$ and $L^2(0,1)$. 

For (nice decaying) 
functions $f:\rone\to\cc$, it 
may be defined via the Fourier transform via $\widehat{\ed f}(\xi)= 
(1+4\pi^2|\xi|^2)^{-1} \hat{f}(\xi)$ or more explicitly, via 
$$
\ed f(x)=e^{-|\cdot|}/2 *f(x) =\f{1}{2}\intl_{-\infty}^\infty  e^{-|x-y|}f(y) dy.
$$
For the case of finite interval, we consider only the case $(0,1)$ 
 for notational convenience. We remark that the results in the 
general case can be recovered by a simple change of variables in the equation. 
Thus for a function $f:(0,1)\to \cc$ given by its Fourier expansion 
$
f(x)= \sum_k a_k e^{2\pi  i k x}$, set  
$$
\ed f (x)= \suml_k \f{a_k}{1+4\pi^2 k^2} e^{2\pi  i k x}
$$

Next, we  verify  that at least formally, the non viscous 
 Camassa-Holm equation \eqref{eq:10} 
 satisfies the conservation law 
$$
\int_I (u^2(t,x) + u_x^2(t,x)) dx=const. 
$$
\subsection{Conservation law for \eqref{eq:10}}
More specifically, let  $I(t)= \int_I (u^2(t,x) + u_x^2(t,x)) dx$. 
If $u$ is a solution, which is 
 sufficiently smooth 
and decaying\footnote{This needs  justification 
in each instance, if  one takes $u$ to a be a solution
of \eqref{eq:10}}, we may  take time  derivative $I'(t)$ to get 
$$
I'(t)= - 2 \int_I (u F(u, u_x) + u_x \p_x [F(u, u_x)]) dx.
$$
\begin{lemma}
\label{le:5} 
Let $u\in C^2(\rone)$, with square integrable second derivative. Then 
$$
\int_I (u F(u, u_x) + u_x \p_x [F(u, u_x)]) dx=0
$$
\end{lemma}
\begin{proof}
This is a simple, although  lengthy  computation. Note that in 
what follows below, all the boundary terms are zero, either 
because $\lim_{|x|\to \infty} u=0$ (in the case $I=\rone$), 
or by the periodic boundary conditions.  We have 
\begin{eqnarray*}
& & \int_I (u F(u, u_x) + u_x \p_x [F(u, u_x)]) dx=  \\
& &= \int_I 
u \left(\f{1}{2}\p_x(u^2)+\p_x\ed [u_x^2/2+u^2]\right) dx + \\
& & + \int_I  u_x \left(\f{1}{2}\p^2_x(u^2)+\p^2_x\ed [u_x^2/2+u^2]\right) dx 
\end{eqnarray*} 
We start with the terms on the second line above. We have 
by integration by parts
\begin{eqnarray*}
& & \f{1}{2} \int_I  u_x \p^2_x(u^2) dx= - \int_I u_{xx} u_x u dx = -\f{1}{2} \int_I \p_x[u_x^2] u dx= 
\f{1}{2} \int_I  u_x^3 dx.
\end{eqnarray*} 
For the second term, use that $\p_x^2\ed= -Id +\ed$ to get 
\begin{eqnarray*}
& &  \int_I  u_x \p^2_x\ed [u_x^2/2+u^2]  dx  = - \int_I u_x[u_x^2/2+u^2]  dx + \\
& &+ 
\int_I u_x\ed[u_x^2/2+u^2]dx = - \f{1}{2} \int_I u_x^3 dx- \int_I u\p_x \ed[u_x^2/2+u^2]dx, 
\end{eqnarray*} 
where we have used that $\int_I u_x u^2 dx=0$. Putting everything together yields the Lemma. 
\end{proof}

\subsection{Littlewood-Paley projections and function spaces} 
Fix  a smooth,  even  
function $\psi\in C_0^\infty(\rone)$, so that $0\leq \psi\leq 1$, $\psi(\xi)=1$, 
 whenever $|\xi|\leq 1$, $\psi$ is decreasing in $(0, \infty)$ 
 and $\psi(\xi)=0$ for all $|\xi|\geq 3/2$. Let also 
$\vp(\xi):=\psi(\xi)-\psi(2\xi)$. Clearly $\vp(\xi)=1$ for all $3/4\leq |\xi|\leq 1$ and 
$\textup{supp} \vp\subset 1/2\leq |\xi|\leq 3/2$. For every integer $k$, define 
the {\it Littlewood-Paley operators}, acting on test functions $f\in \cs(\rn)$ via 
\begin{eqnarray*}
& & \widehat{P_{<2^k} f}(\xi):=\psi(2^{-k} \xi) \hat{f}(\xi),\\
& & 
\widehat{P_{2^k} f}(\xi):=\vp(2^{-k} \xi) \hat{f}(\xi),
\end{eqnarray*}
Clearly the kernels of these operators are given by $2^{kn} \hat{\psi}(2^k \cdot)$ and 
$2^{kn} \hat{\vp}(2^k \cdot)$ respectively and thus commute with differential operators. 
It is also easy to see that since 
$\|2^{kn} \hat{\psi}(2^k \cdot)\|_{L^1}=C\|\hat{\psi}\|_{L^1}$ and similar for the other kernel, {\it $P_{<2^k}, P_{2^k}$ are bounded on $L^p$ spaces for all $1\leq p\leq \infty$ }
with bounds independent of $k$.

The Calder\'on commutator theorem states that the commutator \\
$[P_{2^k}, a] f:=P_{2^k}(a f)- a P_{2^k} f$ acts as a smoothing operator of 
order one\footnote{Similar statement holds for the commutator $[P_{<2^k}, a]$ as well.}. 
More precisely, we shall need a (standard)  estimates of the form  
\begin{eqnarray*}
& & \norm{[P_{2^k}, a] f}{L^r}\leq C 2^{-k} \norm{\nabla a}{L^q}\norm{f}{L^p},\\
& & \norm{[P_{<2^k}, a] f}{L^r}\leq C 2^{-k} \norm{\nabla a}{L^q}\norm{f}{L^p}\\
& & \norm{[P_{2^k}, a]\nabla  f}{L^r}\leq C \norm{\nabla a}{L^q}\norm{f}{L^p},\\
& & \norm{[P_{<2^k}, a]\nabla  f}{L^r}\leq C \norm{\nabla a}{L^q}\norm{f}{L^p},
\end{eqnarray*}
whenever $1\leq r, q, p\leq \infty$ and $1/r=1/p+1/q$.

This whole theory can be developed for the case of finite interval, 
with some notable differences, some of which we discuss below. 

The Littlewood-Paley operators acting on $L^2([0,1])$ are   defined via 
$$
P_{\leq N} f (x) = \suml_{k: |k|\leq N} a_k e^{2\pi i k x},
$$
that is $P_{\leq N}$ truncates the terms in the 
Fourier expansion with frequencies 
$k: |k|> N$.  Clearly $P_{\leq N}$ is a projection 
operator. More generally, we may define for all $0\leq N<M\leq \infty$
$$
P_{N\leq \cdot\leq M}f(x) = \suml_{k: N\leq |k|\leq M} a_k e^{2\pi i k x}.
$$
It is an elementary exercise in orthogonality, that whenever 
$[N_1, M_1]\cap [N_2, M_2]=\emptyset$, then $\int_0^1 
P_{N_1\leq \cdot\leq M_1}f(x)P_{N_2\leq \cdot\leq M_2}g(x)dx=0$. 

For products of three functions, we have the following 
\begin{lemma}
\label{le:80}
Let $f, g, h:[0,1]\to \cc$, with Fourier coefficients $\{f_n\}, \{g_n\}, \{h_n\}$ respectively. 
Then 
$$
\int_0^1 f(x) g(x) h(x) dx=\suml_{m, k\in\cz} f_m g_{-m-k} h_k.
$$
As a consequence, for every $N>>1$,
\begin{equation}
\label{eq:98}
\int_0^1 (P_{>N} f(x)) g(x) (P_{<N/2} h(x)) dx = 
\int_0^1 (P_{>N} f(x)) (P_{>N/2} g(x)) (P_{<N/2} h(x)) dx 
\end{equation}
\end{lemma}
\begin{proof}
The proof is based on the Fourier expansion and the 
fact that $\int_0^1 e^{2\pi i n x} dx =\de_n$. 
More specifically, 
\begin{eqnarray*}
& & \int_0^1 f(x) g(x) h(x) dx = \suml_{m,n, k\in\cz} f_m g_{n} 
h_k \intl_0^1 e^{2\pi i (m+n+k) x} dx = \\
& & =
\suml_{m,n, k\in\cz} f_m g_{n} h_k \de_{m+n+k}= \suml_{m, k\in\cz} f_m g_{-m-k} h_k.
\end{eqnarray*}
For \eqref{eq:98},  observe that if $|m|>N$ and $|k|<N/2$, then $|-m-k|>N/2$. 
\end{proof}
Our next lemma is a well-known Sobolev embedding type result  for the spaces $L^q(0,1)$.  
We state it in the form of the {\it Bernstein inequality}, since this is what we use later on. 
One can also formulate a version in terms of the 
Sobolev spaces defined below. 
\begin{lemma}
\label{le:bern}
Let $N$ be an integer and $f:[0,1]\to \cc$. Then, for every $1\leq p\leq 2 \leq q\leq \infty$, 
$$
\norm{P_{<N} f}{L^q}\leq N^{1/p-1/q} \norm{f}{L^p}.
$$
\end{lemma}
\begin{proof}
First, we establish the lemma for $p=2$, $q=\infty$. Let $f=\sum_n f_n e^{2\pi i n x}$. Then 
$$
\norm{P_{<N} f}{L^\infty}\leq \sum_{n :|n|<N} |f_n|\leq N^{1/2} 
(\sum_{n :|n|<N} |f_n|^2)^{1/2}\leq  N^{1/2} \norm{f}{L^2}. 
$$
Since by Plancherel's theorem $P_{<N}:L^2\to L^2$, it follows that 
$\norm{P_N}{L^q\to L^2}\leq N^{1/2-1/q}$. The rest of the range follows by duality. 
\end{proof}

Introduce some function spaces. Take  
\begin{eqnarray*}
& & 
\dot{H}^s(\rone) = \{f:\rn\to \cc:  (\int_{\rone} |\hat{f}(\xi)|^2|\xi|^{2s}d\xi)^{1/2}<
\infty\}, \\
& & H^s(\rone) = L^2(\rone)\cap \dot{H}^s(\rone), \\
& & \dot{H}^s(0,1) = \{f:(0,1)\to \cc:  
(\suml_{k\in \cz}|a_k|^2|k|^{2s} )^{1/2}<
\infty\}, \\
& & 
H^s((0,1)) = L^2(0,1)\cap \dot{H}^s(0,1).
\end{eqnarray*}
By the Plancherel's theorem $\norm{P_{2^k} f}{\dot{H}^s}\sim 2^{ks} 
\norm{P_{2^k} f}{L^2}$ and  
$\norm{P_{>2^k} f}{\dot{H}^s}\gtrsim  2^{ks} 
\norm{P_{2^k} f}{L^2}$. \\ 
{\bf Remark: }We note that while the Littlewood-Paley operators acting on functions in $L^2(\rn)$ enjoy the Caler\'on commutation estimates, {\it such commutator estimate fails for Littlewood-Paley operators acting on functions in $L^2(I)$}. 

We will also frequently use the fractional differentiation operators of order $s: -\infty<s<\infty$, 
 defined via 
$$
\widehat{|\p|^s f}(\xi):= |\xi|^s \hat{f}(\xi),
$$
in the  case of whole line  and 
$$
|\p|^s (\suml_k a_k e^{2\pi i k x}):= \suml_{k\neq 0}  a_k |k|^s e^{2\pi i k x}. 
$$
in the case $I=(0,1)$.  We would like to point out that $|\p|^s:H^s_0\to L^2_0$ 
is an isometry  and in general 
$$
\norm{|\p|^{s_2} u}{\dot{H}^{s_1}}= \norm{u}{\dot{H}^{s_1+s_2}}.
$$

As a corollary of Lemma \ref{le:bern}, we have that for all $\si>0$, there is 
$C_\si$, so that 
\begin{equation}
\label{eq:sob}
\norm{u}{L^\infty(0,1)}\leq |\int_0^1 u(x) dx|+C_\si\norm{u}{\dot{H}^{1/2+\si}}. 
\end{equation}
\subsection{Kato-Ponce Lemma in the finite interval case}
Recall  the Kato-Ponce product  estimates, that is 
for all $s\geq 0$ and  $1\leq p, q_1, r_1, q_2, r_2 \leq \infty: 1/p=1/q_1+1/r_1=1/q_2+1/r_2$, 
then 
$$
\norm{|\p|^s (f g)}{L^p(\rn)}\leq C_s(\norm{|\p|^s f}{L^{q_1}(\rn)}\norm{g}{L^{r_1}(\rn)}+ 
\norm{f}{L^{r_2}(\rn)}\norm{|\p|^s g}{L^{q_2}(\rn)}).
$$
Unfortunately, we do not know of an analogue of 
such fractional differentiation product estimate for the case of finite interval.  
However, when $s$ is an integer, we have a similar, if somewhat weaker  estimate.  
\begin{lemma}
\label{le:KP}
Let  $s\geq 0$ be an integer 
 and  $1\leq p, q_1, r_1, q_2, r_2 \leq \infty: 1/p=1/q_1+1/r_2=1/q_2+1/r_1$. 
Then for any $X\subset \rn$, 
$$
\norm{\p^s (f g)}{L^p(X)}\leq C_s(\norm{\p^s f}{L^{q_1}(X)}+
\norm{f}{L^{q_2}(X)})( \norm{\p^s g}{L^{r_1}(X)}+\norm{g}{L^{r_2}(X)}).
$$
\end{lemma}
\begin{proof}
Recall the differentiation formula
$$
\p^s (f g)= \suml_{s_1=0}^s \f{s!}{s_1! (s-s_1)!} \p^{s_1} f \p^{s-s_1} g.
$$
and the Young's inequality $ab \leq a^p/p+b^q/q$ for any $1<p,q<\infty: 1/p+1/q=1$. We have 
$$
\norm{\p^s (f g)}{L^p(X)}\leq  2^s \supl_{0\leq s_1\leq s} \norm{\p^{s_1} f \p^{s-s_1} g}{L^p}.
$$
Thus, it will suffice to show that for any $s_1\in[0,s]$, 
$$
\norm{\p^{s_1} f \p^{s-s_1} g}{L^p}\leq C_s (\norm{\p^s f}{L^{q_1}(X)}+
\norm{f}{L^{q_2}(X)})( \norm{\p^s g}{L^{r_1}(X)}+\norm{g}{L^{r_2}(X)}).
$$
Fix $s_1$ and denote $\al=s_1/s\in[0,1]$. If $\al=0$ or $\al=1$, an application of the 
H\"older's inequality gives the result. If $\al\in(0,1)$, then in fact $1/s\leq \al<1-1/s$. \\
Let $\tilde{q}, \tilde{r}$ be
 determined by 
\begin{eqnarray*}
& & 
\tilde{q}^{-1}=\al q_1^{-1}+(1-\al) q_2^{-1},\\
& & 
\tilde{r}^{-1}=(1-\al)r_1^{-1}+\al r_2^{-1}.
\end{eqnarray*} 
Clearly $\tilde{q}^{-1}+\tilde{r}^{-1}=p^{-1}$ and by H\"older's inequality and 
convexity of the norms 
\begin{eqnarray*}
& & 
\norm{\p^{s_1} f \p^{s-s_1} g}{L^p}\leq \norm{\p^{s_1} f }{L^{\tilde{q}}} 
\norm{\p^{s-s_1} g }{L^{\tilde{r}}}\leq 
\norm{\p^{s} f }{L^{q_1}}^\al \norm{f}{L^{q_2}}^{1-\al} 
\norm{\p^{s} g }{L^{r_1}}^{1-\al} \norm{g }{L^{r_2}}^{\al}. 
\end{eqnarray*} 
By Young's inequality, the last expression is bounded by 
$$
C_\al(\norm{\p^{s} f }{L^{q_1}}+\norm{f}{L^{q_2}})(\norm{\p^{s} g }{L^{r_1}}+
\norm{g }{L^{r_2}}),
$$
where $C_\al$ may be taken  $2 \max(\al^{-2}, (1-\al)^{-2})\leq 2 s^2$. 

\end{proof}
\subsection{Attractors}
\label{sec:attractors}
In this section, we offer some basic definitions and 
elementary properties of attractors. 
 
For  an initial value problem for  well-posed evolution equation, 
$$
\frac{d}{dt} u(t)=F(u(t)), \ \ \ u(0)=u_0,
$$ 
defined on a Hilbert space $H$,  consider  the solution semigroup 
$\{S(t)\}_{t \geq 0}$ by $S(t)u_0=u(t)$. $S(t)$ maps $H$ into $H$, satisfies the semigroup
properties $$S(t+s)=S(t) S(s), S(0)=Id$$ and is continuous in the  initial data for 
each $t \geq 0$.

\begin{definition} Let $S(t)$ be a $C_0$ semigroup, acting on a normed space $H$. 
 Then 
\begin{itemize}
\item  
$S(t)$ is called {\it point dissipative} if there is a bounded set $B\subset H$ 
such that for any $u_0 \in H,
S(t)u_0 \in B$ for all sufficiently large $t \geq 0$. That is 
$$
\supl_{u_0\in B}\limsup_{t\to \infty}  \norm{S(t) u_0}{H}<\infty. 
$$
\item 
$S(t)$ is called {\it asymptotically compact} in $H$ if $S(t_n)u_n$ has a 
convergent subsequence for any
{\it bounded sequence} $u_n$ when $t_n \to +\infty$.
\end{itemize}
\end{definition}
\noindent Our next definition gives  a precise
 meaning to the notion of attractor. 
\begin{definition}
$\ca \subset H$ is called a global attractor for the evolution equation if it is compact,
invariant ($S(t) \ca=\ca, \ t \geq 0$) and  attracts every bounded set $X$ 
( $S(t) X \to \ca, \ t \to \infty$).
\end{definition}
A classical result in dynamical systems is that an attractor exists, if 
$S(t)_{t \ge 0}$ is both point dissipative and  asymptotically compact.

Next, we recall the Riesz-Rellich Criteria for precompactness, 
see Theorem XIII.66, p. 248, \cite{SimoniV}). 
\begin{proposition} 
  \label{rrthm}
  Let $S \subseteq L^p(\rn)$ with  $1 \le p <\infty$. Then 
  $S$ is precompact in $L^p(\rn)$ if and only if the following conditions
  are satisfied:
  
  (1) $S$ is bounded in $L^p(\rn)$;
  
  (2) $f \to 0$ in $L^p$ sense at infinity uniformly in $S$, i.e.,
  for any $\epsilon$, there is a bounded set $K \subset \rn$ so that
  for all $f \in S$:
 $ \int_{\rn \backslash K} |f(x)|^p dx \le \epsilon^p; 
 $
 
 (3)  $f(\cdot -y) \to f$ uniformly in $S$ as $y \to 0$, i.e., for any $\epsilon$, there
 is  $\delta$ so that $f\in S$ and $|y|<\delta$ imply that
 $ \int_{\rn} |f(x-y) -f(x)|^p dx \le \epsilon^p.
 $
 \end{proposition}
As shown in \cite{Stanislavova}, \cite{wang}
 (see also Proposition 3 in \cite{Stanislavova1}), we may replace 
the difficult to verify condition $(3)$ in the Riesz-Rellich Criteria above 
by an equivalent condition, which basically says that the ($L^2$ or the $H^1$) 
 mass of the high-frequency component 
has to go uniformly to zero. The exact formulation is 
\begin{proposition} 
\label{prop:3}
Assume that
\begin{itemize}
\item $\supl_n \norm{u_n(t_n,.)}{H^1(\rn)} \leq C $
\item $\limsup\limits_n \norm{u_n(t_n,.)}{H^1(|x|>N)} \to 0 \ {\textup as} \ N \to \infty$
\item $\limsup\limits_n \norm{P_{>N}u_n(t_n,.)}{H^1(\rn)} \to 0 \ {\textup as} \ N \to \infty$
\end{itemize}
Then the sequence  $\{u_n(t_n,.)\}$ is precompact in $H^1(\rn)$.
Same results hold, if one replaces $H^1(\rn)$ by $L^2(\rn)$ everywhere in the statement above. 
\end{proposition}
In the case of finite domains, one has of course the second condition automatically satisfied 
and we have 
\begin{proposition} 
\label{prop:4}
For the  sequence $\{u_n\}\subset H^1(0,1)$, assume 
\begin{itemize}
\item $\supl_n \norm{u_n(t_n,.)}{H^1(0,1)} \leq C $
\item $\limsup\limits_n \norm{P_{>N} u_n(t_n,.)}{H^1(0,1)}\to 0  \ \textup{as}\  N \to \infty$. 
\end{itemize}
then the sequence  $\{u_n(t_n,.)\}$ is precompact in $H^1(0,1)$.
\end{proposition}
We reproduce the short proof of Proposition \ref{prop:4}. 
\begin{proof}
By the Plancherel's theorem, it suffices to show that $b^k=\{a_n^k\}$, $k=1, \ldots$ 
is  
precompact in the weighted space $l^2_s$ if it is uniformly bounded and \\ 
$\lim_{N\to \infty} 
\limsup_k (\sum_{n:|n|>N} |n|^{2s} |a_n^k|^2)^{1/2}=0$. 

By the uniform boundedness of $\{b^k\}$ 
and the reflexivity of $l^2_s$, we have a weak limit $b=\{a_n\}\in l^2_s$ 
of some subsequence 
of $b^k$. Without loss of generality, assume $b^k\to b$ weakly. In particular, for all $n$, 
$a_n^k \to_k a_n$.  
We will show that actually  
$\liml_k\norm{b^k-b}{l^2_s}=0$. 

Fix $\si>0$ and find $N$, so that 
for all,  but finitely many $k$
$$
(\sum_{n:|n|>N} |n|^{2s} |a_n^k|^2)^{1/2}\leq \si/3.
$$
Next, find $N_1$, so that 
$$
(\sum_{n:|n|>N} |n|^{2s} |a_n|^2)^{1/2}\leq \si/3.
$$
Finally, find $k_0$, so that for all $-\max(N, N_1)\leq n\leq \max(N, N_1)$ and for all 
$k>k_0$, we have $|a_n^k-a_n|\leq \si/(10 \max(N, N_1))$. We conclude 
that for all but finitely many $k>k_0$, we have 
$$
\norm{b^k-b}{l^2_s}\leq \si. 
$$
\end{proof}

\section{Global well-posedness for the viscous Camassa-Holm equation}
\label{sec:90}
In this section, we show the global well-posedness for \eqref{eq:2} in both the 
finite interval case and the whole line case. 
The methods are identical in both cases, so we treat it in the same proof. 

As we have mentioned earlier the unbounded operator $A: Au=-\p_x (a(x)u_x)$, 
satisfying \eqref{eq:c1} defines a $C_0$ (and in fact analytic) semigroup, see for example 
\cite{SimonII}, p. 252.   \\
This allows us to reformulate \eqref{eq:2} in an equivalent integral 
equation form\footnote{for smooth and decaying solutions}  
\begin{equation}
\label{eq:15}
u = e^{-t A} u_0-\intl_0^t e^{(s-t)A} F(u, u_x)(s) ds.
\end{equation}
Our first step then will be to show a local well-posedness result. 
\subsection{Local well-posedness for \eqref{eq:2}} 
\label{sec:lwp}
Regarding the simpler equation 
\eqref{eq:11}, we have taken the classical 
 approach for the heat equation outlined in \cite{McOwen}. 
We will use the following lemma, which is a compilation of 
Theorem 3  (p. 298-300) and the discussion in Section 11.2.b, \cite{McOwen}.
\begin{lemma}
\label{le:17}
Suppose $S(t)=e^{-t L}$ is a $C_0$-semigroup acting on both $L^2(I)$ and 
$\dot{H}^1(I)$. 
Assume also   
\begin{eqnarray}
\label{eq:c10}
& & \norm{S(t) g}{\dot{H}^1(I)}\leq C t^{-1/2} \norm{g}{L^2}.
\end{eqnarray}
For  the integral equation 
$$
u(t)= S(t) u_0 +\intl_0^t S(t-s) F(u)(s) ds,
$$
there exists time $T>0$ depending only on $\norm{u_0}{H^1}$, 
such that the integral equation  has an unique local solution 
$u\in C([0,T], H^1)$ provided  
\begin{equation}
\label{eq:163}
\left\{\begin{array}{l}
\norm{F(u)-F(v)}{L^2}\leq M_R\norm{u-v}{H^1} \\
\textup{whenever}\qq \norm{u}{H^1},\norm{v}{H^1}\leq R. 
\end{array}
\right.
\end{equation}
\end{lemma} 
We first show the proof of Lemma \ref{le:17} and then verify 
\eqref{eq:c10}  for 
the semigroup $S(t)=e^{-tA}$ and \eqref{eq:163} for the Camassa-Holm nonlinearity $F(u, u_x)$. '
\begin{proof}(Lemma \ref{le:17})
We set a fixed point argument for the integral equation at hand. Set 
$X^R_T=\{u\in C([0,T), H^1(I)), \sup_{0<t<T}\norm{u(t, \cdot}{H^1}\leq R \}$ and the map 
$$
\La u(t, \cdot) = S(t) u_0 +\intl_0^t S(t-s) F(u)(s) ds.
$$
We need to show that for appropriate $R=R(\norm{u_0}{H^1})$ and  $T=T(R)$, 
$\La:X^R_T\to X^R_T$ is a contraction. Take $R=10 \norm{u_0}{H^1}$. To see  
$\La:X^R_T\to X^R_T$, we have by \eqref{eq:c10} and  \eqref{eq:163} (applied for the case $v=0$),  
\begin{eqnarray*}
& &\norm{\La u(t, \cdot)}{H^1}\leq \norm{S(t)u_0}{H^1}+
\|\intl_0^t e^{(s-t)L} F(u)(s) ds\|_{L^2}+  
C \|\intl_0^t e^{(s-t)L} F(u)(s) ds\|_{\dot{H}^1} \\
& & \leq C\norm{u_0}{H^1}+
\intl_0^t \norm{F(u)(s, \cdot)}{L^2}ds + \intl_0^t \f{\norm{F(u)(s, \cdot)}{L^2}}{\sqrt{t-s}} ds\leq \\
& & \leq  C\norm{u_0}{H^1}+ C M(\norm{u}{H^1}) (t+\sqrt{t})\supl_{0<t<T} \norm{u}{H^1}. 
\end{eqnarray*}
Clearly choosing $T=T(R)$ small enough,  $0<t<T$  and $\sup_{0<t<T} \norm{u}{H^1}\leq R$ 
 will guarantee that the right hand side is less than $R$. One verifies similarly the 
contraction property of $\La:X^R_T\to X^R_T$, by using the full strength of \eqref{eq:163}. 
\end{proof}

First, we verify that $e^{-tA}$ is a semigroup on $\dot{H}^1$. 
Observe that $\norm{u}{\dot{H}^1}\sim \norm{A^{1/2} u}{L^2}$. Indeed, 
$$
\norm{A^{1/2} u}{L^2(I)}^2 = \dpr{Au}{u}=\int_I  a(x) u_x^2 dx\sim \norm{u_x}{L^2}^2,
$$
by \eqref{eq:c1}. Then 
$$
\norm{e^{-t A} f}{\dot{H}^1}\sim \norm{A^{1/2} e^{-t A} f}{L^2}= 
\norm{ e^{-t A} A^{1/2}  f}{L^2}\leq C\norm{A^{1/2}  f}{L^2}\sim \norm{f}{\dot{H}^1}. 
$$

The estimate \eqref{eq:c10} is a  standard property of analytic 
semigroups, see Corollary 1 and Corollary 2, \cite{SimonII}, p. 252. We choose to 
 deduce it as a 
 simple consequence of the functional calculus for the 
self adjoint operator $A$. 

We have 
$\norm{e^{-tA} g}{\dot{H}^1}\sim \norm{A^{1/2} e^{-tA} g}{L^2} = t^{-1/2} 
\norm{f(tA) g}{L^2}$, where 
$f(y)=e^{-y}y^{1/2}$ is a well-defined bounded function on the spectrum of $A$. 
It follows that 
$$
\norm{e^{-tA} g}{\dot{H}^1}\leq Ct^{-1/2} \norm{f}{L^\infty(0,\infty)} \norm{g}{L^2}
\leq Ct^{-1/2} \norm{g}{L^2}, 
$$ 
which is \eqref{eq:c10}. 

It remains to establish  \eqref{eq:163} for the Camassa-Holm nonlinearity $F$. We actually 
prove a little more general statement. 
\begin{lemma}
\label{le:20}
Let $F$ be  the nonlinearity  for the Camassa-Holm equation, as defined earlier. 
Then for all nonnegative  integers $s$, we have 
\begin{equation}
\label{eq:164}
\norm{F(u)-F(v)}{H^s}\leq M (\norm{u}{H^{s+1}}+\norm{v}{H^{s+1}}) \norm{u-v}{H^{s+1}}. 
\end{equation} 
\end{lemma}
\begin{proof}
We have by Lemma \ref{le:KP} and the Sobolev embedding $L^\infty(I)\hookrightarrow H^{1/2+}(I)
\hookrightarrow H^{s+1}(I)$,  
\begin{eqnarray*}
& & \norm{\p_x(u^2)-\p_x(v^2)}{\dot{H}^s}\sim \norm{|\p_x^{s+1}[(u-v)(u+v)]}{L^2}\leq \\
& &\leq 
C( \norm{|\p_x^{s+1} (u-v)}{L^2}+\norm{(u-v)}{L^\infty})(\norm{u+v}{L^\infty}+  
\norm{|\p_x^{s+1}(u+v)}{L^2})\leq \\
& & \leq M \norm{u-v}{H^{s+1}} (\norm{u}{H^{s+1}}+\norm{v}{H^{s+1}}). 
\end{eqnarray*}
For the second term in $F$, consider first $s\geq 1$. We have 
\begin{eqnarray*}
& & \norm{\p_x\ed(u_x^2-v_x^2)}{\dot{H}^s}\leq C 
\norm{|\p_x^{s-1}[ (u_x-v_x) (u_x+v_x)]}{L^2} \leq \\
& & \leq 
C(\norm{|\p_x^{s}(u-v)}{L^\infty}+ \norm{u_x-v_x}{L^2})
(\norm{u}{H^1}+\norm{v}{H^1}+
\norm{\p_x^s u}{L^\infty}+ \norm{\p_x^s v}{L^\infty})\leq \\
& & \leq C \norm{u-v}{H^{s+1/2+}}( \norm{u}{H^{s+1/2+}}+\norm{v}{H^{s+1/2+}})\leq  
C \norm{u-v}{H^{s+1}}( \norm{u}{H^{s+1}}+\norm{v}{H^{s+1}}). 
\end{eqnarray*} 
When $s=0$, use either  
$$
\p_x\ed f(x)=\f{1}{2} \intl sgn(x-y) e^{-|x-y|} f(y) dy,
$$
or 
$$
\p_x\ed f(x)=\suml_n \f{2\pi i n}{1+4\pi^2 n^2} f_n e^{2\pi i  n x},
$$
to conclude that 
$\p_x\ed:L^1(I)\to L^2(I)$. It follows that   
$$
\norm{\p_x\ed(u_x^2-v_x^2)}{L^2(I)}
\lesssim \norm{(u_x-v_x)(u_x+v_x)}{L^1}\lesssim \norm{u-v}{H^1}(\norm{u}{H^1}+\norm{v}{H^1}). 
$$
For the third term in $F$, we easily estimate 
\begin{eqnarray*}
& & \norm{\p_x\ed(u^2-v^2)}{\dot{H}^s}\leq C \norm{u-v}{H^{\max(s-1,0)}}(\norm{u}{L^\infty}+ 
\norm{v}{L^\infty})\leq \\
& &\leq C \norm{u-v}{H^{s+1}}(\norm{u}{H^{s+1}}+\norm{u}{H^{s+1}}).
\end{eqnarray*}
\end{proof}
Note that  one can represent $F(u)=\La(u,u)$, where $\La(u,v)$ is the  bilinear form 
$$
\La(u,v)= \f{1}{2}\p_x(u v)+\p_x\ed [u_x v_x/2+u v].
$$
It is easy to see that one can show (with the same exact proof)  for every integer $s\geq 0$ 
$$
\norm{\La(\vp,\psi)}{H^s}\leq C\norm{\vp}{H^{s+1}}\norm{\psi}{H^{s+1}}.
$$
A bilinear interpolation between the estimates above, (which are  valid for all integers), yields 
the corresponding estimates for non integer values of $s$ as well. Setting $\vp=\psi=u$, we obtain
\begin{corollary}
\label{cor:1}
Let $s\geq 0$ and $F$ be the Camassa-Holm nonlinearity. Then 
$$
\norm{F(u)}{H^s}\leq M \norm{u}{H^{s+1}}^2. 
$$
\end{corollary}
\subsection{$H^2$ smoothness of the local solutions}  
\label{sec:smoothness}
In this section, we show the $H^2$ smoothness of the local $H^1$ solution  
constructed above. Beside the obvious importance of having this extra smoothness 
information, this will enable us (see Section \ref{sec:global} below) to iterate the 
local solution to a global one by utilizing the conservation (or rather dissipation) 
of the $H^1$ energy. We have 
\begin{proposition}
\label{prop:1}
Let $u$ be the $H^1$  solution to \eqref{eq:15}, with  life span $T$. 
Then there exists a  constant $C_\ve$, so that for all $0<t<T$, $u\in C((0,t), H^2(I))$ and 
as a result 
$$
\norm{u(t, \cdot)}{H^2(I)}\leq \f{C_\ve}{\sqrt{t}} \norm{u_0}{H^1}+ C_\ve t^{1/4} 
\supl_{0<s<t} 
\norm{u(s, \cdot)}{H^1}^{3}+ C_\ve\norm{u(t, \cdot)}{H^1}. 
$$
\end{proposition}
\begin{proof}
The argument required for the proof is to rerun again the fixed point method, this time in the smoother space $H^2(I)$. However, this amounts to showing $H^2$ {\it a priori} estimates for the solution, which is what we concentrate on. 

 Apply $A$ to \eqref{eq:15}. This is justified, since the right hand side of \eqref{eq:15} 
is in the domain of $A$ by the 
semigroup properties of $e^{-t A}$. 
We have 
\begin{eqnarray*}
\norm{Au}{L^2}\leq \norm{e^{-tA} A u_0}{L^2}+C \intl_0^t \norm{e^{(s-t)A} A F(u)(s)}{L^2(I)} ds
\end{eqnarray*}
But  
$$
\norm{e^{-tA} A u_0}{L^2}= \norm{e^{-tA} A^{1/2} (A^{1/2} u_0)}{L^2}\leq C t^{-1/2} 
\norm{A^{1/2} u_0}{L^2}\sim C t^{-1/2} \norm{u_0}{H^1}. 
$$
On the other  hand, by the properties of the  functional calculus for $A$
\begin{equation}
\label{eq:901}
\norm{e^{-z A} A F}{L^2}=|z|^{-1} \norm{f(A) F}{L^2}\leq C |z|^{-1} (\supl_{y>0} |e^{-y}y| )
\norm{F}{L^2}\leq C |z|^{-1} \norm{F}{L^2},
\end{equation}
for all $z>0$, while 
\begin{equation}
\label{eq:902}
\norm{e^{-z A} A F}{L^2}\leq C\norm{A F}{L^2}\leq C_\ve \norm{F}{H^2}. 
\end{equation}
The last inequality can be checked easily as follows 
\begin{eqnarray*}
& & \norm{A u}{L^2}^2=\int (a u_{xx}+a' u_x)^2 dx= \int (a^2 u_{xx}^2 - a a'' u_x^2) dx\leq \\
& &\leq 
\norm{a}{L^\infty}^2 \norm{u_{xx}}{L^2}^2+\norm{a}{L^\infty}\norm{a''}{L^\infty}
\norm{u_{x}}{L^2}^2\lesssim \norm{u}{H^2}^2.
\end{eqnarray*}
A complex interpolation  between \eqref{eq:901} and \eqref{eq:902}  yields 
$$
\norm{e^{-z A} A F}{L^2}\leq C |z|^{-7/8} \norm{F}{H^{1/4}}. 
$$
Plugging this estimate back in the integral term yields 
\begin{eqnarray*}
& & \intl_0^t \norm{e^{(s-t)A} A F(u)(s)}{L^2(I)} ds\leq C \intl_0^t 
\f{\norm{F(s, \cdot)}{H^{1/4}}}{(t-s)^{7/8}}ds\leq C t^{1/8} \sup_{0<s<t}
\norm{F(s, \cdot)}{H^{1/4}} .
\end{eqnarray*}
According to Corollary \ref{cor:1},  $\norm{F(s, \cdot)}{H^{1/4}}\leq C \norm{u}{H^{5/4}}^2$. By the 
Gagliardo-Nirenberg inequality, 
$\norm{u}{H^{5/4}}\leq \norm{u}{H^2}^{1/4}\norm{u}{H^1}^{3/4}$. 

Putting everything together 
\begin{eqnarray*}
\norm{Au(t, \cdot)}{L^2} &\leq & \f{C}{\sqrt{t}}\norm{u_0}{H^1}+ C t^{1/8} \supl_{0<s<t}
\norm{u(s, \cdot)}{H^2}^{1/2}\norm{u(s, \cdot)}{H^1}^{3/2}\leq \\
& &\leq 
\f{C}{\sqrt{t}}\norm{u_0}{H^1}+ C_\si t^{1/4}  \supl_{0<s<t} 
\norm{u(s, \cdot)}{H^1}^{3}+ \si \supl_{0<s<t}  \norm{u(s, \cdot)}{H^2}.
\end{eqnarray*}
for any $\si>0$ and some $C_\si$. 

Observe now,
$\norm{Au}{L^2}^2\geq \ve^2 \norm{u}{H^2}^2/2-C \norm{u}{H^1}^2$. Indeed, 
\begin{eqnarray*}
& & \norm{Au}{L^2}^2= \int (a^2 u_{xx}^2 +(a')^2 u_x^2 +2 a a' u_x u_{xx}) dx\geq \\
& & \geq \int \f{a^2}{2}
 u_{xx}^2 dx - \int (a')^2 u_x^2 dx \geq \ve^2 \norm{u}{H^2}^2/2-C \norm{u}{H^1}^2.
\end{eqnarray*}
Let $G(t)=\supl_{0<s\leq t}  \norm{u(s, \cdot)}{H^2}$. Taking into account the last inequality 
provides 
$$
G(t)\leq \f{C_\ve}{\sqrt{t}}\norm{u_0}{H^1}+ C_{\si,\ve} t^{1/4}  \supl_{0<s<t} 
\norm{u(s, \cdot)}{H^1}^{3}+ C\norm{u(t, \cdot)}{H^1}+ C_\ve\si G(t).
$$
Choosing 
appropriately small $\si:C_\ve\si<1/2 $, allows us to hide the last term and  as a result 
$$
 \norm{u(t, \cdot)}{H^2}\leq G(t)\leq \f{C_{\ve}}{\sqrt{t}}\norm{u_0}{H^1}+ C_\ve t^{1/4} 
\supl_{0<s<t} 
\norm{u(s, \cdot)}{H^1}^{3}+ C_\ve\norm{u(t, \cdot)}{H^1}. 
$$
\end{proof}
\noindent {\bf Remark} The above argument can be extended (with no additional 
smoothness or otherwise assumptions on $A$) to show that $u\in \cap_{m=0}^\infty D(A^m)$ with the 
corresponding estimates (away from the zero) for $\norm{A^m u}{L^2}$ 
as in Proposition \ref{prop:1}. This is  
the usual regularity result that one expects for parabolic equations. 

\subsection{Global well-posedness for \eqref{eq:2}} 
\label{sec:global}
Our approach to global well-posedness for the  parabolic problem 
\eqref{eq:2} is to iterate the 
local well-posedness result to a global one.

We will  show that for the local $H^1$ solution, produced in Section \ref{sec:lwp}, 
one has the estimate 
\begin{equation}
\label{eq:21}
\norm{u(t, \cdot)}{H^1}\leq I(0) e^{C t} + C_\ve(e^{C t}-1) 
\supl_{0\leq s\leq t} \norm{g(s, \cdot)}{L^2}^2.
\end{equation}
for every $0<t<T$, where $T$ is its lifespan. 

Assuming \eqref{eq:21}, let us prove that the solution is global. 
Fix $u_0\in H^1(I)$ and define for every (sufficiently large) 
integer $n$
$$
T_n=\sup\{ t:H^1 \textup{solution is defined in } 
(0,t)\ \&\  \supl_{0<t_1<t} \norm{u(t_1, \cdot)}{H^1}<n\},
$$
and $T^*=\limsup_n T_n$. 

If $T^*=\infty$, there is nothing to prove, 
the solution is global. If $T^*<\infty$, it must be  that 
$\limsup_{t\to T^*}\norm{u(t, \cdot)}{H^1}=\infty$. On the other hand,  take any 
 sequence $t_n\to T^*$. By \eqref{eq:21}, 
$$
\limsup_{n\to \infty} \norm{u(t_n, 
\cdot)}{H^1}\leq I(0) e^{C T^*} + C_\ve(e^{C T^*}-1) 
\supl_{0\leq s\leq T^*} \norm{g(s, \cdot)}{L^2}^2 <\infty,
$$
a contradiction. This implies the solutions produced in Section \ref{sec:lwp} are global ones. 
Therefore, it remains to show \eqref{eq:21}. 
\subsubsection{Local boundedness of $t\to \norm{u(t, \cdot)}{H^1}$}
In view of the $H^2$ smoothness, established in Proposition \ref{prop:1}, 
this follows in a  standard 
way from Lemma \ref{le:5}. To this end, let 
$$
I(t)=\int_I (u^2(t,x)+u_x^2(t,x)) dx.
$$ 
and differentiate in time. Then one may use the equation (because of the $H^2$ smoothness) to get  
\begin{eqnarray*}
& & I'(t)= 2 \int_I (u u_t + u_x (u_t)_x) dx = -2 \int (u F(u, u_x)+u_x \p_x F(u, u_x)) dx + \\
& & + 2\int_I 
u \p_x (a(x) u_x)dx +2 \int_I 
u_x \p^2_x (a(x) u_x)dx + 2\int_I u g(t,x) dx+ 2\int_I u_x g_x(t,x) dx 
\end{eqnarray*}
Note that by Lemma \ref{le:5}, $\int (u F(u, u_x)+u_x \p_x F(u, u_x)) dx=0$. For the next term, 
clearly 
$$
\int_I 
u \p_x (a(x) u_x)dx= -\int a(x) u_x^2 dx\leq 0
$$
Next, consider the term $ \int_I 
u_x \p^2_x (a(x) u_x)dx$. We have 
\begin{eqnarray*}
& & \int_I 
u_x \p^2_x (a(x) u_x)dx =- \int \p_x(a u_x) u_{xx} dx= - \int a(x) u_{xx}^2 +
\f{1}{2}\int a''(x) u_x^2 dx\leq \\
& &\leq  -\ve \norm{u_{xx}}{L^2}^2+ \norm{a''}{L^\infty} \norm{u_x}{L^2}^2\leq 
-\ve \norm{u_{xx}}{L^2}^2+ C \norm{u_x}{L^2}^2 
\end{eqnarray*}
We have used here $a(x)\geq \ve$ and  $a\in C^2(I)$. 

Finally, we have 
\begin{eqnarray*}
& &|\int_I u g(t,x) dx+ \int_I u_x g_x(t,x) dx |\leq C(\norm{u}{L^2}+\norm{u_{xx}}{L^2})
\norm{g(t, \cdot)}{L^2}\leq \\
& & \leq  \ve \norm{u_{xx}}{L^2}^2/2+ \norm{u}{L^2}^2+
C_\ve \norm{g(t, \cdot)}{L^2}^2. 
\end{eqnarray*}
Altogether, 
\begin{eqnarray*}
& &I'(t)\leq -\ve \norm{u_{xx}}{L^2}^2/2+C(\norm{u(t, \cdot)}{L^2}^2+\norm{u_{x}(t, \cdot)}{L^2}^2)+
C_\ve\norm{g(t, \cdot)}{L^2}^2\leq \\
& & \leq C I(t)+C_\ve\norm{g(t, \cdot)}{L^2}^2
\end{eqnarray*}
Rewrite this as 
\begin{eqnarray*}
& & \f{d}{d t} (e^{-C t} I(t))\leq C_\ve e^{-C t} \norm{g(t, \cdot)}{L^2}^2,
\end{eqnarray*}
whence upon integration we get 
$$
I(t)\leq I(0) e^{C t} + C_\ve(e^{C t}-1) \supl_{0\leq s\leq t} \norm{g(s, \cdot)}{L^2}^2. 
$$
which  is  \eqref{eq:21}.

\section{Global attractors for the viscous Camassa-Holm: The finite interval case}
\label{sec:attractor_1}
In this section, we prove Theorem \ref{theo:2}. As we have discussed in 
Section \ref{sec:attractors} and more specifically Proposition \ref{prop:4}, 
we will need to verify that for any $t_n\to \infty$ and for any $B>0$ and 
any sequence of initial data $\{u_n\}\subset H^1(0,1)$ with 
$\supl_n \norm{u_n}{H^1}\leq B$, we have 
\begin{eqnarray}
\label{eq:70}
& & \sup_{u_0\in H^1_0} \limsup_{t\to \infty}\norm{S(t)u_0}{H^1}\leq C(g, \ve), \\
\label{eq:71}
& & \sup_n \norm{S(t_n)u_n}{H^1}\leq C(B; g, \ve), \\
\label{eq:72}
& & \lim_N \limsup_n\norm{P_{>N} S(t_n)u_n}{H^1}=0.
\end{eqnarray}
This section is devoted to showing \eqref{eq:70}, \eqref{eq:71}. The estimate 
 \eqref{eq:72} is somewhat more complicated and it will postponed 
until Section \ref{sec:vanish}. 
In the end, we will 
show the asymptotic smoothing effect, that is the fact that the attractor lies in a smoother space. 
\subsection{Point dissipativeness: Proof of \eqref{eq:70}}
Fix $u_0$ with $\norm{u_0}{H^1}\leq B$. Consider the solution to \eqref{eq:2} with initial data $u_0$, 
$u(t, \cdot)=S(t)u_0$. We have already shown the local boundedness of 
$t\to \norm{u(t, \cdot)}{H^1}$ ( i.e. is  \eqref{eq:21}), 
 which we now improve. Note that the 
extra conditions $\int_0^1 g(x)dx=\int_0^1 u_0(x)dx=0$ are  crucial in our argument. 

Recall $I(t)=\int_I (u^2(t,x)+u_x^2(t,x)) dx$. We need to reexamine our estimates above for $I'(t)$, in order to use to our advantage the smallness of $\norm{a'}{L^\infty}$.  We have as before 
\begin{eqnarray*}
& & I'(t)= 2 \int_I (u u_t + u_x (u_t)_x) dx =  - 2\int_I 
a(x) (u_x)^2 dx -2 \int_I 
u_{xx} \p_x (a(x) u_x)dx + \\ 
& & +2\int_I u g(t,x) dx- 2\int_I u_{xx} g(t,x) dx \leq -2\ve(\norm{u_x}{L^2}^2 +
\norm{u_{xx}}{L^2}^2)+ \\
& &+ 2 \norm{a'}{L^\infty} \norm{u_x}{L^2} \norm{u_{xx}}{L^2}+ 
(\norm{u}{L^2}+\norm{u_{xx}}{L^2})\norm{g}{L^2}. 
\end{eqnarray*}
Since $\norm{a'}{L^\infty}\leq \ve$, it is easy to see that the term 
$2 \norm{a'}{L^\infty} \norm{u_x}{L^2} \norm{u_{xx}}{L^2}$ gets absorbed  
 by $\ve(\norm{u_x}{L^2}^2 +
\norm{u_{xx}}{L^2}^2)$ and we get 
\begin{equation}
\label{eq:22}
I'(t)\leq -\ve( \norm{u_x}{L^2}^2+\norm{u_{xx}}{L^2}^2) + 
(\norm{u}{L^2}+\norm{u_{xx}}{L^2})\norm{g}{L^2}
\end{equation}

Note that by the conservation law 
$\p_t \int_0^1 u(t,x) dt= \int_0^1 g(x) dt=0$ and $\int_0^1 u_0(x)=0$, 
we have $\int_0^1 u(t,x)dx=0$ for all $t$. Let $u(t,x)=\sum_{n\neq 0} a_n(t) e^{2\pi i n x}$. 
It follows that 
$$
\norm{u(t, \cdot)}{L^2}=\left(\suml_{|n|\geq 1} |a_n|^2\right)^{1/2}\leq 
\left(\suml_{|n|\geq 1} |n|^2 |a_n|^2\right)^{1/2}\leq C  \|u_x(t, \cdot)\|_{L^2}.
$$
Use $\norm{u(t, \cdot)}{L^2}\leq  \|u_x(t, \cdot)\|_{L^2}\leq  \|u_{xx}(t, \cdot)\|_{L^2}$ 
and the 
  Cauchy-Schwartz's inequality in \eqref{eq:22} to get  
\begin{eqnarray*}
& & I'(t)\leq  -\ve (\norm{u_x}{L^2}^2+\norm{u_{xx}}{L^2}^2) + 
(\norm{u_x}{L^2}+\norm{u_{xx}}{L^2})\norm{g}{L^2}\leq  \\
& & \leq -\ve( \norm{u}{L^2}^2+\norm{u_{xx}}{L^2}^2)/2 + 
C \norm{g}{L^2}^2/\ve  \leq -\ve I(t)/2 + 
C\norm{g}{L^2}^2/\ve 
\end{eqnarray*}
We now finish with a Gronwall type argument, namely we rewrite the inequality above as 
$$
\f{d}{dt} (I(t) e^{t \ve /2}) \leq C e^{t \ve /2} 
 \norm{g}{L^2}^2/\ve,
$$
which after integration in time yields 
\begin{equation}
\label{eq:23}
I(t)\leq I(0) e^{-  \ve t/2} + C \norm{g}{L^2}^2/\ve^2. 
\end{equation}
It follows that 
$$
\limsup_{t\to \infty} I(t) \leq C \norm{g}{L^2}^2/\ve^2,
$$
which is the point dissipativeness of $S(t)$. 

\subsection{Uniform boundedness: Proof of \eqref{eq:71} } The uniform boundedness 
in fact follows from \eqref{eq:23} as well. Indeed,  denote 
$I_n(t)=\norm{S(t)u_n}{L^2}^2+\norm{S(t)u_n}{\dot{H}^1}^2$. Clearly $I_n(0)=\norm{u_n}{H^1}^2\leq B^2$. 
We have by \eqref{eq:23}, 
$$
I_n(t_n)\leq I_n(0) e^{- \ve t_n/2} + C \norm{g}{L^2}^2/\ve^2 \leq 
B^2+C \norm{g}{L^2}^2/\ve^2. 
$$
\section{Uniform vanishing: Proof of \eqref{eq:72}}
\label{sec:vanish}
Fix a real number $B$. Let the  initial data be 
$u_0: \norm{u_0}{H^1}\leq B$, 
with a corresponding solution $u$. 
We know from the results of the previous sections that 
such solutions exist globally and belong to the class $C((t_1, t_2), H^2)$ 
for every $0<t_1<t_2<\infty$. 

Let $k$ be a (large) positive integer and denote  
$$
I_{>k}(t)=\intl_0^1 ((P_{>2^k} u)^2 + (P_{>2^k} u_x)^2 dx.
$$
This is the high-frequency portion of the energy, which we are trying 
to show is small as $N\to \infty$, uniformly in $\norm{u_0}{H^1}$. 
We use energy estimate reminiscent of the 
 estimate for $I(t)$. 

After taking  time derivative, use the equation \eqref{eq:2}  and $P_{>2^k}^2=P_{>2^k}$. We get 
\begin{eqnarray*}
& & I_{>k}'(t)=2\intl_0^1  ( P_{>2^k} u P_{>2^k} u_t  + P_{>2^k} u_x P_{>2^k} u_{tx} dx= \\
& &= 2 \intl_0^1 P_{>2^k} u F(u, u_x)+ P_{>2^k} u_x \p_x F(u, u_x) dx+ \\
& &+ \intl_0^1 
P_{>2^k} u \p_x( a(x)  u_x)dx+ P_{>2^k} u_x \p^2_x( a(x)  u_x)dx+ \\
& & +
\intl_0^1 
(P_{>2^k} u g + P_{>2^k} u_x g_x dx)=:N+V+F
\end{eqnarray*}
There are three sort of terms arising in the energy estimate. 
We start with those  arising from the viscosity.   
\subsection{Viscosity terms}
\label{sec:vis}
Write 
\begin{eqnarray*}
& & \f{V}{2}= \intl_0^1 
(P_{>2^k} u) \p_x( a(x)  u_x)dx+ (P_{>2^k} u_x) \p^2_x( a(x)  u_x)dx= \\
& & = - \intl_0^1 
\left[(P_{>2^k} u_x) a(x)  u_x dx+ (P_{>2^k} u_{xx})a(x) u_{xx}+ (P_{>2^k} u_{xx})a'(x) u_{x} \right]dx.
\end{eqnarray*}
We estimate the first and the third term  by H\"older's inequality and the uniform boundedness 
$$
|\intl_0^1 
(P_{>2^k} u_x) a(x)  u_x dx|\leq \norm{a}{L^\infty} \norm{u_x}{L^2}^2\leq C(B; g, \ve).
$$
Also, by H\"older and Cauchy-Schwartz
\begin{eqnarray*}
& & |\intl_0^1  (P_{>2^k} u_{xx})a'(x) u_{x} dx|\leq \norm{a'}{L^\infty} \norm{u_x}{L^2} 
\norm{P_{>2^k} u_{xx}}{L^2}\leq \\
 & & \leq \f{\ve}{100} \norm{P_{>2^k} u_{xx}}{L^2}^2+ 
\f{C}{\ve} \norm{a'}{L^\infty}^2 \norm{u_x}{L^2}^2=
\f{\ve}{100} \norm{P_{>2^k} u_{xx}}{L^2}^2+C(B; g, \ve).
\end{eqnarray*}
We need more delicate estimates for the second term $\int (P_{>2^k} u_{xx})a(x) u_{xx} dx$. 
The difficulties here lie with the fact that the commutators 
$[P_{>N},a]$ {\it are not smoothing operators }, when considered on $L^2[0,1]$,  
(in contrast with $L^2(\rone)$). 

Write $u_{xx}= P_{>2^k} u_{xx}+P_{\leq 2^k}u_{xx}$  to get 
$$
\int (P_{>2^k} u_{xx})a(x) u_{xx} dx= \int (P_{>2^k} u_{xx})^2 a(x) dx+
\int (P_{>2^k} u_{xx})a(x) P_{\leq 2^k}u_{xx} dx.
$$
Clearly, $\int (P_{>2^k} u_{xx})^2 a(x) dx\geq \ve \norm{P_{>2^k} u_{xx}}{L^2}^2$, while we will show 
\begin{equation}
\label{eq:50}
\begin{array}{c}
|\int (P_{>2^k} u_{xx})a(x) P_{\leq 2^k}u_{xx} dx|\leq  \\
\leq 
C 2^k (\norm{a'}{L^\infty} \norm{P_{>2^{k-1}}u_{x}}{L^2} 
\norm{P_{>2^k} u_{xx}}{L^2}+ 
\norm{a_{>{2^{k-1}}}}{L^\infty}\norm{P_{>2^k} u_{xx}}{L^2} \norm{u_x}{L^2}). 
\end{array}
\end{equation}
To that end, write 
\begin{eqnarray*}
& &
\int (P_{>2^k} u_{xx})a(x) P_{\leq 2^k}u_{xx} dx= \int (P_{>2^k} u_{xx})a(x) P_{2^{k-1}< \cdot 
\leq 2^k}u_{xx} dx+ \\
& &+\int (P_{>2^k} u_{xx})a(x) P_{\leq  2^{k-1}}u_{xx} dx. 
\end{eqnarray*}
For the first term, 
use that $a(x)=a(0)+\intl_0^x a'(y) dy$ and by orthogonality \\ 
$\int (P_{>2^k} u_{xx})a(0) P_{2^{k-1}< \cdot 
\leq 2^k}u_{xx} dx=0$. We get 
\begin{eqnarray*}
& & |\int (P_{>2^k} u_{xx})a(x) P_{2^{k-1}\leq \cdot 
\leq 2^k}u_{xx} dx| = |\int (P_{>2^k} u_{xx})(\int_0^x a'(y)dy)  P_{2^{k-1}<\cdot 
\leq 2^k}u_{xx} dx|\leq \\
& & \leq \norm{a'}{L^\infty} \norm{P_{>2^k} u_{xx}}{L^2} 
\norm{P_{2^{k-1}< \cdot 
\leq 2^k} u_{xx}}{L^2} \leq 2^k \norm{a'}{L^\infty} \norm{P_{2^{k-1}< \cdot 
\leq 2^k} u_{x}}{L^2}
\norm{P_{>2^k} u_{xx}}{L^2} \leq \\
& & \leq 2^k \norm{a'}{L^\infty} \norm{P_{>2^{k-1}} u_{x}}{L^2}
\norm{P_{>2^k} u_{xx}}{L^2} . 
\end{eqnarray*}
For the second term, use Lemma \ref{le:80}, 
more specifically \eqref{eq:98}. We have 
\begin{eqnarray*}
& & |\int (P_{>2^k} u_{xx})a(x) P_{\leq  2^{k-1}}u_{xx} dx| =  
|\int (P_{>2^k} u_{xx})(P_{>2^{k-1}} a(x)) P_{\leq  2^{k-1}}u_{xx} dx|\leq \\
& &\leq  \norm{P_{>2^k} u_{xx}}{L^2} \norm{P_{>2^{k-1}} a}{L^\infty} 
\norm{P_{\leq  2^{k-1}}u_{xx}}{L^2}\leq 
C 2^{k} \norm{a_{>2^{k-1}}}{L^\infty} \norm{P_{>2^k} u_{xx}}{L^2} \norm{P_{\leq  2^{k-1}} u_x}{L^2} \\
& & \leq C 2^{k} \norm{a_{>2^{k-1}}}{L^\infty} \norm{P_{>2^k} u_{xx}}{L^2} \norm{u_x}{L^2}.
\end{eqnarray*}
This establishes \eqref{eq:50}. 

Put together all terms
 that arise from the viscosity and use  the uniform boundedness 
 \eqref{eq:71} and  the Cauchy-Schwartz inequality  $a b\leq \ve a^2 +(4\ve)^{-1}b^2$
to obtain 
\begin{eqnarray*}
& & 
V\leq  - \f{2\ve}{3}   \norm{P_{>2^k} u_{xx}}{L^2}^2+ 
C(B; g,\ve, \de) + \f{2^{2k} \norm{a'}{L^\infty}^2}{\ve} 
\norm{P_{>2^{k-1}}u_{x}}{L^2}^2+  \\
& & +
 2^{2k} 
\norm{a_{>2^{k-1}}}{L^\infty}^2 C(B; g,\ve)\leq  - \f{2\ve}{3}   \norm{P_{>2^k} u_{xx}}{L^2}^2+ 2^{2k} \de^2 \ve 
\norm{P_{>2^{k-1}}u_{x}}{L^2}^2+  C(B; g,\ve) + \\
& &+
2^{2k} 
\norm{a_{>2^{k-1}}}{L^\infty}^2 C(B; g,\ve).
\end{eqnarray*}
The last inequality holds due to $\norm{a'}{L^\infty}\leq \de \ve$. 
\subsection{Nonlinearity terms}
\label{sec:nonl} For the nonlinearity terms, we have several easy terms, that we take care of first. 
Namely, according to Lemma \ref{le:20} (see \eqref{eq:164} with $s=0$)
\begin{eqnarray*}
& & |\int (P_{>2^k} u) F(u,x_x) dx|\leq \norm{ P_{>2^k} u}{L^2}\norm{F(u, u_x)}{L^2}\leq 
C \norm{u}{H^1}^3\leq C(B; g,\ve). 
\end{eqnarray*}
Also, by H\"older's inequality and the Sobolev embedding  
\eqref{eq:sob}
\begin{eqnarray*}
& & |\int (P_{>2^k} u_x) \p_x^2 (u^2) dx|= |\int (P_{>2^k} u_{xx}) \p_x (u^2) dx| \leq 
\norm{P_{>2^k} u_{xx}}{L^2} \norm{u_x}{L^2}\norm{u}{L^\infty}\leq \\
& &\leq C \norm{P_{>2^k} u_{xx}}{L^2} \norm{u}{H^1}^2\leq \f{\ve}{100} \norm{P_{>2^k} u_{xx}}{L^2}^2+ 
C(B; g,\ve).
\end{eqnarray*}

Finally,  
\begin{eqnarray*}
& & |\int P_{>2^k} u_x \p_x^2 \ed(u_x^2/2+u^2) dx|\leq C \norm{P_{>2^k} u_x}{L^\infty} \norm{u}{H^1}^2.
\end{eqnarray*}
However,  by Lemma \ref{le:bern}
\begin{eqnarray*}
& & \norm{P_{>2^k} u_x}{L^\infty}\leq \suml_{l\geq k} 
\norm{P_{2^l<\cdot\leq 2^{l+1}} u_x}{L^\infty} \leq 
\suml_{l\geq k} 2^{l/2} 
\norm{P_{2^l<\cdot\leq 2^{l+1}} u_x}{L^2} \sim \\
& & \sim 
\suml_{l\geq k} 2^{-l/2} 
\norm{P_{2^l<\cdot\leq 2^{l+1}} u_{xx}}{L^2}\leq C 2^{-k/2} \norm{P_{>2^k} u_{xx}}{L^2}\leq 
C \norm{P_{>2^k} u_{xx}}{L^2}, 
\end{eqnarray*}
implying that 
$$
|\int P_{>2^k} u_x \p^2 \ed(u_x^2/2+u^2) dx|\leq \f{\ve}{100} \norm{P_{>2^k} u_{xx}}{L^2}^2 
+  C(B; g,\ve). 
$$

\subsection{Forcing terms}
\label{sec:forc} The forcing terms are easy  to control. 
\begin{eqnarray*}
& & 
|\int P_{>2^k} u g + P_{>2^k} u_x  g_x dx|= |\int P_{>2^k} u g -  P_{>2^k} u_{xx}  g dx|\leq \\
& & \leq (\norm{P_{>2^k} u}{L^2}+\norm{P_{>2^k} u_{xx}}{L^2}) \norm{g}{L^2}\leq \\
& & \leq 
\f{\ve}{100} (\norm{P_{>2^k} u}{L^2}^2+\norm{P_{>2^k} u_{xx}}{L^2}^2 )+ 
\f{C}{\ve} \norm{g}{L^2}^2\leq \\
& &\leq \f{\ve}{50} \norm{P_{>2^k} u_{xx}}{L^2}^2 + C(B; g, \ve).
\end{eqnarray*}

\subsection{Conclusion of the argument for uniform vanishing of the high frequencies} 
Put together all the estimates for viscosity terms, 
forcing terms and nonlinearity terms. We obtain 
\begin{eqnarray*}
& &  I_{>k}'(t)\leq -\f{\ve}{2}\norm{P_{>2^k} u_{xx}}{L^2}^2+ C 2^{2k} \de^2 \ve 
 \norm{P_{>2^{k-1}} u_{x}}{L^2}^2+ \\
& &+
2^{2k}\norm{a_{>2^{k-1}}}{L^\infty}^2  C(B;  g, \ve)+ 
C(B;  g, \ve).
\end{eqnarray*}
Note first that $\norm{P_{>2^k} u_{xx}}{L^2}\geq 2^k \norm{P_{>2^k} u_{x}}{L^2}\geq c
2^k \sqrt{I_{>k}(t)}$. 

Next, we estimate the term $\norm{a_{>2^{k-1}}}{L^\infty}$. Let 
$a(x)=\suml_l a_l e^{2\pi i l x}$. Then  
\begin{eqnarray*}
& & \norm{a_{>2^{k-1}}}{L^\infty}\leq \suml_{l>2^k} 
|a_l| \leq 
C 2^{-k} \suml_{l>k} |l ||a_l| \leq C 2^{-k} (\suml_{l>k} 
|a_l|^2 |l|^{4})^{1/2}(\suml_{l>1} l^{-2})^{1/2} \leq \\
& & \leq C 2^{-k} \norm{a''}{L^2(I)}\leq C 2^{-k} \norm{a''}{L^\infty}. 
\end{eqnarray*}
We plug in this estimate to get 
\begin{equation}
\label{eq:100}
I_{>k}'(t)+ \f{2^{2k} \ve}{4}I_{>k}(t)\leq C 2^{2k} \de^2 \ve 
 \norm{P_{>2^{k-1}} u_{x}}{L^2}^2+ C(B;  g, \ve)
\end{equation}
Notice that as before, we can rewrite \eqref{eq:100} as 
$$
\f{d}{d t} (I_{>k}(t) e^{2^{2k} \ve t/4}) \leq C 2^{2k} \de^2 \ve e^{2^{2k}  \ve t/4}
 \norm{P_{>2^{k-1}} u_{x}}{L^2}^2+ e^{2^{2k}  \ve t/4} C(B;  g, \ve),
$$
which after time integration yields 
\begin{equation}
\label{eq:104}
\begin{array}{l}
I_{>k}(t)\leq I_{>k}(0)e^{-2^{2k} \ve t/4}+ C \de^2 \supl_{0\leq s\leq t}  
\norm{P_{>2^{k-1}} u_{x}(s, \cdot)}{L^2}^2+ 2^{-2k} C(B;  g, \ve)\leq \\
\leq I_{>k}(0)e^{-2^{2k} \ve t/4}+ C \de^2 \supl_{0\leq s\leq t}  I_{>k-1}(s)
+ 2^{-2k} C(B;  g, \ve). 
\end{array}
\end{equation}
Informally, it should be that $I_{>k-1}\sim I_{>k}$, and since $\de^2<<1$, 
we may ignore the middle term and get the desired uniform vanishing. However, $I_{>k-1}\geq I_{>k}$ and 
we may not perform this operation. \\
To go around this difficulty, introduce  
$$
I^n_{>k}(t)= \int 
((u^n_{>2^k}(t, \cdot))^2+(\p_x u^n_{>2^k}(t, \cdot))^2)dx, 
$$ 
where $\{u^n\}\subset H^1$, with $\sup_n \norm{u^n}{H^1}\leq B$. Note that by the uniform 
boundedness \eqref{eq:71}, we have 
$$
\sup_{n, k, t} I^n_{>k}(t)\leq \int 
((u^n(t, \cdot))^2+(\p_x u^n(t, \cdot))^2)dx\leq C(B; g, \ve).
$$
Let also 
$h^n_k(t)= \sup_{0\leq s\leq t} I^n_{>k}(s)$. Recast  \eqref{eq:104} for each $n$  as 
\begin{equation}
\label{eq:117}
h^n_k(t)\leq h^n_k(0)e^{-2^{2k} \ve t/4} + C\de^2 h^n_{k-1}(t)+ C 2^{-2k} C(B;  g, \ve)
\end{equation}
We will need $\de$ so small, that $C\de^2\leq 1/8$. Denote also $h_k=\limsup_{n\to \infty} 
h^n_k(t_n)$ for some fixed sequence $t_n\to \infty$. Thus, we have 
$$
h_k\leq h_{k-1}/8+ 2^{-2k} C(B;  g, \ve)
$$
Iterating this inequality, we obtain 
\begin{eqnarray*}
& & 
h_k\leq h_{k-1}/8+ 2^{-2k} C(B;  g, \ve)\leq 
8^{-2} h_{k-2}+ (2^{-2k}+2^{-2k-1}) C(B;  g, \ve)\leq  \ldots \\
& & \leq  2^{-2k} C(B;  g, \ve)+ 
2^{-3k} h_0\leq (2^{-2k}+ 2^{-3k}) C(B;  g, \ve), 
\end{eqnarray*}
since by \eqref{eq:71},  $h_0\leq C(B;  g, \ve)$. It follows that 
$$
\lim_{k\to \infty} \limsup_{n\to\infty} \norm{P_{>2^k} u^n(t_n, \cdot)}{H^1}=\lim_{k\to \infty}h_k =0, 
$$
which is \eqref{eq:72}. Moreover, we have that the attractor 
(whose existence is now established) is actually a 
{\it bounded subset} of $H^{2-\si}$ for all $\si>0$. 

Indeed, since 
 every element of the attractor is of the form $u(\cdot)=\lim_n u^n(t_n, \cdot)$, we have 
by the last estimate 
$$
\sup_k 2^{2k} \norm{P_{>2^k} u}{L^2}\leq C(B, g, \ve),
$$
or $u\in B^2_{2, \infty}$. 
Of course, this implies 
$$
\norm{u(\cdot)}{H^s}^2\sim \suml_{k\geq 1} 2^{2k(s-1)}
\norm{P_{\sim 2^k} u(\cdot)}{H^1}^2 \leq \suml_{k\geq 1}  2^{2k(s-1)} 2^{-2k} C(B;  g, \ve, \de)<
C(B;  g, \ve)
$$
if $s<2$. 

\section{Attractors for the viscous Camassa-Holm equation on the whole line}
\label{sec:conclusions}
In this section, we indicate the main steps for the Proof of Theorem \ref{theo:9}. 
Since most of the arguments are quite similar to 
those already presented for the case of finite interval, 
we will  frequently refer to the previous 
sections.

To start with, let us point out that Theorem \ref{theo:1}, which applies to the (undamped) viscous Camassa-Holm equation \eqref{eq:2} applies as stated to \eqref{e:1} as well. 
The reader may reproduce the arguments from Section \ref{sec:90} easily, but we point out 
that the energy estimates in fact work better in  the presence of the 
damping factor $\mu u$, see the discussion regarding the proof of \eqref{e:2} below. 

To establish the asymptotic compactness 
of the dynamical system $S(t)$  associated with \eqref{e:1}, 
we resort to Proposition \ref{prop:3}, just as we have 
 used the similar Proposition \ref{prop:4} for the case of finite interval. 

Therefore, fix a sequence of times $\{t_n\}$ and 
$\{u_n\}\subset H^1(\rone)$, which is uniformly bounded, say 
$\sup_n \norm{u_n}{H^1}\leq B$. It remains to show 
\begin{eqnarray}
\label{e:2}
& & \sup_{f\in H^1} \limsup_{t\to \infty}\norm{S(t)f}{H^1}\leq C(g, \mu, \ve) \\
\label{e:3}
& & \supl_n \norm{S(t_n)u_n}{H^1}\leq C(B, g, \ve, \mu) \\
\label{e:4}
& &  \lim_{N\to \infty} \limsup_n \norm{P_{>N} S(t_n) u_n}{H^1}=0 \\
\label{e:5}
& & \lim_{N\to \infty} \limsup_n \norm{S(t_n) u_n}{H^1(|x|>N)}=0.
\end{eqnarray}
Note that \eqref{e:2} is the point dissipativeness of $S(t)$, while \eqref{e:3},\eqref{e:4}, 
 \eqref{e:5} guarantee the asymptotic compactness of $S(t)$, according 
to Proposition \ref{prop:3}. 
\subsection{Proof of \eqref{e:2}}
Denote 
$
I(t)= \int_{\rone} (u^2 +u_x^2) dx
$
and compute 
\begin{eqnarray*}
& & I'(t)= 2\int u u_t + u_x u_{x t} dx = 2\int u(-F(u, u_x) +\p_x(a u_x)-\mu u+g) dx - \\
& & - 2
\int u_{xx} ( -F(u, u_x) +\p_x(a u_x)-\mu u+g) dx= \\
& &=-2 \int (u F(u, u_x) +u_x \p_x F(u, u_x)) dx-2 \int a(x) (u_x^2+ u_{xx}^2) dx - \\
& & - 2\int a'(x) u_{xx} u_x dx+2\int (u-u_{xx}) g -2 \mu \int (u^2+u_x^2) dx 
\end{eqnarray*}
We split now our considerations, depending on the assumptions on $a$.\\
{\bf Estimate with the assumption $\norm{a'}{L^\infty}<<\ve$.}\\ 
By Lemma \ref{le:5}, $\int (u F(u, u_x) +u_x \p_x F(u, u_x)) dx=0$ and we 
estimate the rest by H\"older's inequality 
\begin{eqnarray*}
& & 
I'(t)\leq -2\ve \int (u_x^2+ u_{xx}^2) dx+2 \norm{a'}{L^\infty} \norm{u_x}{L^2} 
\norm{u_{xx}}{L^2}+\\
& & +2 \norm{g}{L^2}(\norm{u}{L^2}+\norm{u_{xx}}{L^2}) - 2 \mu \int (u^2+u_x^2) dx 
\end{eqnarray*}
By  the smallness of $\norm{a'}{L^\infty}$, we conclude \\
 $\norm{a'}{L^\infty} \norm{u_x}{L^2} 
\norm{u_{xx}}{L^2}\leq \ve( \norm{u_x}{L^2}^2+\norm{u_{xx}}{L^2}^2)/2$. 
On the other hand, by Young's inequality 
$$
\norm{g}{L^2}(\norm{u}{L^2}+\norm{u_{xx}}{L^2})\leq \mu \norm{u}{L^2}^2/2+ 
\ve\norm{u_{xx}}{L^2}^2/4+\f{C}{\min(\mu, \ve)}\norm{g}{L^2}^2.
$$ 
Altogether, 
\begin{equation}
\label{e:89}
I'(t)\leq -\f{\ve}{2} \int (u_x^2+u_{xx}^2 )dx-\mu \int u^2 dx+ \f{C}{\min(\mu, \ve)}\norm{g}{L^2}^2.
\end{equation}
We show that \eqref{e:89} follows by assuming $a''(x)\leq 2 a(x)$. \\
{\bf Estimate with the assumption $2a''(x)\leq 2 a(x)$.}\\
We perform one more integration by parts in the expression for $I'(t)$ to get 
\begin{eqnarray*}
& & I'(t)=  - 2 \int a(x) (u_x^2+ u_{xx}^2) dx +\int a''(x) u^2_x dx+2\int (u-u_{xx}) g - \\
& & - 2 \mu \int (u^2+u_x^2) dx \leq -2 \int a(x) u_{xx}^2 -2\mu \int (u_x^2 +u^2) dx 
+2\int (u-u_{xx}) g dx \leq \\
& & \leq -\min(\ve, \mu)  \int (u^2 + u_x^2) dx  + \f{C}{\min(\mu, \ve)} \norm{g}{L^2}^2.
\end{eqnarray*}

Thus, under either the smallness assumption $\norm{a'}{L^\infty}<<\ve$ or under $a''(x)\leq 2 a(x)$, we have 
$$
I'(t)+\f{\min(\ve, \mu)}{2}  I(t)\leq \f{C}{\min(\mu, \ve)}\norm{g}{L^2}^2,
$$
which by Gronwall's inequality implies 
$$
I(t)\leq e^{- \min(\ve, \mu)/2 t} I(0)+ \f{C}{\min(\mu, \ve)^2}\norm{g}{L^2}^2= 
e^{- \min(\ve, \mu)/2 t} \norm{f}{H^1}^2 + \f{C}{\min(\mu, \ve)^2}\norm{g}{L^2}^2.
$$
Taking limit $t\to \infty$ establishes \eqref{e:2}.
\subsection{Proof of \eqref{e:3}}
Uniform boundedness of the orbits follows from the last estimate as follows. Denote
 $I_n(t)=\norm{u_n(t, \cdot)}{H^1}^2$. We have 
$$
I_n(t)\leq e^{- \min(\ve, \mu)/2 t} \norm{u_n(0, \cdot)}{H^1}^2 + 
\f{C}{\min(\mu, \ve)^2}\norm{g}{L^2}^2\leq B^2+ \f{C}{\min(\mu, \ve)^2}\norm{g}{L^2}^2,
$$
where $B=\sup_n \norm{u_n(0)}{H^1}$. 
\subsection{Proof of \eqref{e:4}}
The proof of \eqref{e:4} largely follows the 
argument for the similar estimate \eqref{eq:72}. Set 
$$
I_{>2^k} (t)= \int_{\rone} (u_{>2^k})^2+(\p_x u_{>2^k})^2 dx
$$
and compute as in Section \ref{sec:vanish} 
\begin{eqnarray*}
& & I_{>2^k}'(t)= 2 \intl P_{>2^k}^2 u F(u, u_x)+ P_{>2^k}^2 u_x \p_x F(u, u_x) dx+ \\
& &+ 2 \int  
P_{>2^k} u \p_x P_{>2^k}( a(x)  u_x)dx+ P_{>2^k} u_x \p^2_x P_{>2^k}( a(x)  u_x)dx+ \\
& & +
2 \int 
(P_{>2^k}^2 u g + P_{>2^k}^2 u_x g_x dx)- 2 \mu \int ((P_{>2^k}^2 u)^2 + (P_{>2^k}^2 u_x)^2 dx).
\end{eqnarray*}

The estimates for the terms arising from the nonlinearity 
 work just in the case of finite interval. Again the damping terms can be ignored,
 because they give  rise to terms with negative signs. 
 
In short, the estimates that we need can be summarized in 
\begin{eqnarray*}
& & |\intl P_{>2^k}^2 u F(u, u_x)+ P_{>2^k}^2 u_x \p_x F(u, u_x) dx| \leq 
C \norm{P_{>2^k} u_{xx}}{L^2}\norm{u}{H^1}^2\leq \\
& & \leq  \ve 
\norm{P_{>2^k} u_{xx}}{L^2}^2+ \f{C}{\ve} \norm{u}{H^1}^4.
\end{eqnarray*}
Similarly, the estimates for the terms arising form the forcing $g$ are estimated by 
\begin{eqnarray*}
& & |\int 
(P_{>2^k}^2 u g + P_{>2^k}^2 u_x g_x )dx|\leq  \ve 
\norm{P_{>2^k} u_{xx}}{L^2}^2/100 +  \f{C}{\ve} \norm{g}{L^2}^2.
\end{eqnarray*}
Finally, the viscosity terms are in fact better behaved than the corresponding terms for the finite interval case, but one has to proceed in a 
slightly different fashion, due to the 
technical inconvenience that $P_{>2^k}$ are not involutions, i.e. 
$P_{>2^k}^2\neq P_{>2^k}$. 

We have 
\begin{eqnarray*}
& & V = - 2 \int  
(P_{>2^k} u_x) P_{>2^k}( a(x)  u_x)dx  -2 \int (P_{>2^k} 
u_{xx}) P_{>2^k}( a(x)  u_x)dx = \\
& & = -2 \int  
(P_{>2^k} u_x)  a(x) (P_{>2^k} u_x)dx-
2 \int (P_{>2^k} u_x) [P_{>2^k},a] u_x dx - \\
& & - 2 \int (P_{>2^k} u_{xx}) P_{>2^k} (a' u_x) dx - 2\int P_{>2^k} u_{xx} 
a(x) P_{>2^k} u_{xx} dx - \\
& & - 2 \int P_{>2^k} u_{xx} [P_{>2^k},a] u_{xx} dx\leq -2\ve 
\int (P_{>2^k} u_{x})^2  + (P_{>2^k} u_{xx})^2 dx + \\
& & + 2 \int |(P_{>2^k} u_x) [P_{>2^k},a] u_x| dx+ 2 
\int |(P_{>2^k} u_{xx}) P_{>2^k} (a' u_x)| dx+ \\
& & + 
2 \int |P_{>2^k} u_{xx} [P_{>2^k},a] u_{xx} |dx. 
\end{eqnarray*}
Note that 
\begin{eqnarray*}
& & 
\int (P_{>2^k} u_{x})^2  + (P_{>2^k} u_{xx})^2 dx\geq 
\int (P_{>2^k} u_{xx})^2 dx \geq c 2^{2k} 
\int (P_{>2^k} u_{x})^2 dx\sim  \\
& & \sim 2^{2k} \int (P_{>2^k} u)^2+ (P_{>2^k} u_{x})^2 dx= 2^{2k} 
I_{>2^k}(t). 
\end{eqnarray*}

By the Calder\'on commutator estimates 
\begin{eqnarray*}
& & \int |(P_{>2^k} u_x) [P_{>2^k},a] u_x| dx\leq C 2^{-k} 
\norm{P_{>2^k} u_x}{L^2} \norm{a'}{L^\infty} \norm{u_x}{L^2}\leq 
C \norm{u_x}{L^2}^2  \norm{a'}{L^\infty}, \\
& & \int |(P_{>2^k} u_{xx}) P_{>2^k} (a' u_x)| dx\leq 
\norm{P_{>2^k} u_{xx}}{L^2}\norm{a'}{L^\infty} \norm{u_x}{L^2}\leq \\
& & \leq \ve 
\norm{P_{>2^k} u_{xx}}{L^2}/100+ \frac{C}{\ve}
\norm{a'}{L^\infty}^2 \norm{u_x}{L^2}^2,   \\
& & \int |P_{>2^k} u_{xx} [P_{>2^k},a] u_{xx} |dx\leq C 
\norm{P_{>2^k} u_{xx}}{L^2} \norm{a'}{L^\infty} \norm{u_x}{L^2}\leq \\
& & \leq 
\ve 
\norm{P_{>2^k} u_{xx}}{L^2}/100+ \frac{C}{\ve}
\norm{a'}{L^\infty}^2 \norm{u_x}{L^2}^2. 
\end{eqnarray*}
Altogether, the various terms in $I_{>2^k}'$ are estimated by 
$$
I_{>2^k}'(t)\leq -\ve 2^{2k} I_{>2^k}(t)+ \frac{C}{\ve}
\norm{a'}{L^\infty}^2 \norm{u_x}{L^2}^2+
\f{C}{\ve}(\norm{g}{L^2}^2+\norm{u}{H^1}^4).
$$
By the uniform boundedness(i.e.  \eqref{e:3} ) and the 
Gronwall's inequality, we deduce 
$$
I_{>2^k}(t)\leq I_{>2^k}(0) e^{-\ve 2^{2k} t}+ 2^{-2k} 
\f{C}{\ve^2}(\norm{g}{L^2}^2+\supl_{0\leq s\leq t}\norm{u(s, \cdot)}{H^1}^4+\norm{a'}{L^\infty}^2\supl_{0\leq s\leq t}\norm{u(s, \cdot)}{H^1}^2 ).
$$
It follows that 
$$
\limsup_n \norm{P_{>2^k} S(t_n) u_n}{H^1}\leq 2^{-k} C(B, g, \ve)
$$
and therefore $\lim_{k\to \infty} \limsup_n \norm{P_{>2^k} S(t_n) 
u_n}{H^1}=0$, thus establishing \eqref{e:4}. 

Note that since $\limsup_n \norm{P_{>2^k} S(t_n) u_n}{H^1}
\lesssim 2^{-k}$, it follows with the same argument as before 
that the attractor 
$\ca\subset H^{2-\si}(\rone)$ for every $\si>0$. 

\subsection{Proof of \eqref{e:5}}
Our last goal is to establish the uniform smallness of the $H^1$ 
energy functional away from large balls. Set 
$$
J_{>N}(t)=\int (u^2(t,x)+u_x^2(t,x))(1-\psi(x/N)) dx.
$$
Compute the  derivative 
\begin{eqnarray*}
& & J_{>N}'(t)= 2\int (u u_t+u_x u_{xt})(1-\psi(x/N)) dx = \\
& &=-2\int (u F(u, u_x)+u_x \p_x F(u, u_x)) (1-\psi(x/N)) dx- \\ 
& &- 2\mu\int 
(u^2+u_x^2)(1-\psi(x/N)) dx+\\
& & +2 \int (u \p_x(a u_x)+ u_x \p_x^2 (a u_x)) (1-\psi(x/N)) dx. 
\end{eqnarray*}
The first term has already been handled  in our previous paper, 
\cite{Stanislavova}. According to  Lemma 5,  \cite{Stanislavova}
the estimate 
 is\footnote{This is actually not so hard to justify. Observe that by 
the conservation law 
$\int (u F(u, u_x)+u_x \p_x F(u, u_x))dx=0$, 
all the integration by parts in \eqref{e:8} that does not hit the 
term $(1-\psi(x/N))$ equates to zero. Therefore, the only terms that survive are those with $N^{-1}\psi'(x/N)$ in them. 
Observe that there are no $u_{xx}$ in those either, whence \eqref{e:8}.}
\begin{equation}
\label{e:8}
|\int (u F(u, u_x)+u_x \p_x F(u, u_x)) (1-\psi(x/N)) 
dx|\leq \f{C}{N}\norm{u(t, \cdot)}{H^1}^3. 
\end{equation}
Next, integration by parts yields 
\begin{eqnarray*}
& & \int (u \p_x(a u_x)+ u_x \p_x^2 (a u_x)) (1-\psi(x/N)) dx = \\
& & = 
-\int a u_x^2 (1-\psi(x/N)) dx+ N^{-1} 
\int a u u_x \psi'(x/N) dx-\\ 
& & - \int u_{xx} \p_x (a u_x) (1-\psi(x/N)) dx+ 
N^{-1} \int u_x \p_x(a u_x)\psi'(x/N) dx 
\end{eqnarray*}
The terms with the factor $N^{-1}$ are  ``good'' terms. 

For the first term , we estimate right away
$$
 N^{-1} 
|\int a u u_x \psi'(x/N) dx|\leq C N^{-1}\norm{a}{L^\infty} \norm{u}{H^1}^2. 
$$
For the second term containing $N^{-1}$, we have 
\begin{eqnarray*}
& & N^{-1} \int u_x \p_x(a u_x)\psi'(x/N) dx = N^{-1} \int a'(x) u_x^2 
\psi'(x/N) dx - \\
& & - \f{1}{2N} \int u_x^2 \p_x( a' \psi'(x/N)) dx\leq 
\f{C}{N} \norm{u_x}{L^2}^2(\norm{a'}{L^\infty}+\norm{a''}{L^\infty}),
\end{eqnarray*}
for some absolute constant $C$. \\
Finally, we have to estimate the term $- \int u_{xx} 
\p_x (a u_x) (1-\psi(x/N)) dx$. As before, we need to use 
either the smallness of $\norm{a'}{L^\infty}$ or $a''(x)\leq 2 a(x)$. \\
{\bf Estimate under the assumption $a''(x)\leq 2 a(x)$.}\\
We have 
\begin{eqnarray*}
& & - \int u_{xx} \p_x (a u_x) (1-\psi(x/N)) dx= \\
& & =
-\int a u_{xx}^2 (1-\psi(x/N)) dx- \int u_{xx} a' u_x (1-\psi(x/N)) dx\leq  \\
& & \leq \f{1}{2} \int u_x^2 \p_x 
(a'(1-\psi(x/N)) ) dx= \\
& &=
\f{1}{2} \int u_x^2 a''(x) (1-\psi(x/N))  dx-\f{1}{2 N} 
\int a''(x) u_x^2 \psi'(x/N) dx\leq \\
& & \leq \f{1}{2} \int u_x^2 a''(x) 
(1-\psi(x/N))  dx+ \f{1}{2N}\norm{u_x}{L^2}^2 \norm{a''}{L^\infty}. 
\end{eqnarray*}
All in all, we get 
\begin{eqnarray*}
& &J_{>N}'(t)\leq 
\f{C}{N}(\norm{u(t, \cdot)}{H^1}^3+ 
(\norm{a}{L^\infty} +\norm{a'}{L^\infty}+\norm{a''}{L^\infty})
\norm{u(t, \cdot)}{H^1}^2)+\\
& &+ \int u_x^2 a''(x) 
(1-\psi(x/N))  dx-2\int a(x) u_x^2 (1-\psi(x/N))  dx -\\
& &- 2\mu \int (u^2 +u_x^2) (1-\psi(x/N))
dx. 
\end{eqnarray*}
We now use the condition 
$a''(x)\leq 2 a(x)$,  to conclude  that the middle term is non-positive whence
$$
J_{>N}'(t)\leq 
\f{C}{N}(\norm{u(t, \cdot)}{H^1}^3+ 
(\norm{a}{L^\infty} +\norm{a'}{L^\infty}+\norm{a''}{L^\infty})
\norm{u(t, \cdot)}{H^1}^2) - 2\mu J_{>N}(t).
$$
By the uniform bounds on 
$\norm{u(t, \cdot)}{H^1}$ , (i.e. \eqref{e:3}) and the previous 
considerations, it follows that 
\begin{equation}
\label{eq:end}
J_{>N}'(t)+ \mu J_{>N}(t)\leq \f{C(B, g, \ve, \mu)}{N}.
\end{equation}
We will show that \eqref{eq:end} holds, by assuming appropriate smallness of 
$\norm{a'}{L^\infty}$. \\
{\bf Estimate under the assumption $\norm{a'}{L^\infty}<<\ve$.}\\
We have 
\begin{eqnarray*}
& & - \int u_{xx} \p_x (a u_x) (1-\psi(x/N)) dx\leq  \\
& & \leq \int |u_{xx}| |a'(x)||u_x| (1-\psi(x/N)) dx- \int u_{xx}^2 a(x)(1-\psi(x/N)) dx \leq \\
& &\leq \norm{a'}{L^\infty} (\int u_{xx}^2 (1-\psi(x/N)) dx )^{1/2}
(\int u_x^2(1-\psi(x/N)) dx )^{1/2}-\\
& &-
\ve \int u_{xx}^2 (1-\psi(x/N)) dx\leq 
2\f{\norm{a'}{L^\infty}^2}{\ve} \int u_x^2(1-\psi(x/N)) dx\leq 2\de^2 \ve J_{>N}(t).
\end{eqnarray*}
Taking $\de$ so small that $\de^2\ve<\mu$ ensures that $2\de^2 \ve J_{>N}(t)$ is 
subsumed by  \\ 
$-2\mu \int u^2+u_x^2)(1-\psi(x/N)) dx$ and therefore, 
we arrive at  \eqref{eq:end} again. \\ 
The Gronwall's inequality  applied to \eqref{eq:end} yields 
$$
J_{>N}(t)\leq e^{-\mu t} J_{>N}(0)+ 
\f{C(B, g, \ve, \mu)}{N}. 
$$
Thus $\limsup_{t_n\to\infty} J_{>N}(t_n)\leq N^{-1} C(B, g, \ve, \mu)$, whence 
$$
\lim_{N\to\infty} \limsup_{t_n\to\infty} J_{>N}(t_n)=0.
$$

\end{document}